\documentclass[a4paper,11pt, draft]{amsart}

\usepackage{amsmath, amsthm, amssymb, amsfonts}

\usepackage{xy}

\newtheorem{theo}{Theorem}[section]

\newtheorem{df}[theo]{Definition}

\newtheorem{lem}[theo]{Lemma}

\newtheorem{prop}[theo]{Proposition}

\newcommand{\R}{\mathbf{R}}
\newcommand{\C}{\mathbf{C}}
\newcommand{\Z}{\mathbf{Z}}
\newcommand{\F}{\mathbf{F}}

\newcommand{\N}{\mathbf{N}}

\newcommand{\la}{\lambda}
\newcommand{\s}{\sigma}
\renewcommand{\a}{\alpha}
\renewcommand{\b}{\beta}
\newcommand{\g}{\gamma}

\renewcommand{\k}{\mathcal{K}}

\newcommand{\G}{\Gamma}
\renewcommand{\t}{\theta}
\newcommand{\bh}{B(\mathcal{H})}
\newcommand{\h}{\mathcal{H}}

\renewcommand{\and}{\mbox{ and }}

\newcommand{\fonction}[4]
{\begin{array}{rcl}
#1 & \to & #2 \\
#3 & \mapsto & #4 \\
\end{array}}

\newcommand{\Core}{\operatorname{Core}}
\newcommand{\Ad}{\operatorname{Ad}}
\newcommand{\Sd}{\operatorname{Sd}}

\newcommand{\id}{\operatorname{id}}
\newcommand{\Aut}{\operatorname{Aut}}

\newcommand{\Out}{\operatorname{Out}}

\newcommand{\centre}{\operatorname{Z}}
\newcommand{\Cent}{\operatorname{Cent}}
\newcommand{\Id}{\operatorname{Id}}
\newcommand{\Mat}{\operatorname{Mat}}
\newcommand{\Proj}{\operatorname{Proj}}
\newcommand{\Inn}{\operatorname{Inn}}

\begin{document}

\title[A new construction of factors of type ${\rm III_1}$]{A New CONSTRUCTION OF FACTORS OF TYPE ${\rm III_1}$}

\author{Cyril Houdayer}

\address{Institut de Math\'ematiques de Jussieu \\
                 Alg\`ebres d'Op\' erateurs et Repr\'esentations \\
                 175 rue du Chevaleret \\
                 75013 Paris \\
                 France}

\email{houdayer@math.jussieu.fr}

\urladdr{http://www.institut.math.jussieu.fr/\textasciitilde houdayer/}

\subjclass[2000]{46L10}

\keywords{Factors of type ${\rm III}$; Full factors; Continuous decomposition.}

\newcommand{\indice}[1]{{\vphantom{#1}}_{\mathit \alpha}{#1}}
\newcommand{\indicebis}[1]{{\vphantom{#1}}_{\mathit \alpha'}{#1}}

\begin{abstract}
We give in this paper a new construction of factors of type ${\rm III_1}$. Under certain assumptions, we can, thanks to a result by Popa, give a complete classification for this family of factors. Although these factors are never full, we can nevertheless, in many cases, compute Connes' $\tau$ invariant. We obtain a new example of an uncountable family of pairwise non-isomorphic factors of type ${\rm III_1}$ with the same $\tau$ invariant.
\end{abstract}

\maketitle

\section*{Introduction}

Let $N$ be a type ${\rm II_{\infty}}$ factor endowed with a trace-scaling one-parameter automorphism group denoted by $(\a_s)_{s \in \R^*_+}$. Let $\Gamma$ be a countable, dense subgroup of $\R^*_+$. The crossed product of $N$ by $\G$ under the action $\a$, denoted by $N \rtimes_{\a} \G$, is a factor of type ${\rm III_1}$; this factor has almost-periodic states and its invariant $\Sd$ is included in $\G$ (see \cite{connes74} for further details).
We want to generalize this construction to the case of a virtual subgroup of $\R$: for $(X, \mu)$ a measure space with $\mu$ a finite or infinite measure, and  $(\g_t)_{t \in \R}$
a measure-preserving, free, ergodic action of $\R$ on $(X, \mu)$ by automorphisms, we shall give the right definition of the crossed product of $N$ by the virtual subgroup $(L^{\infty}(X), \g)$.

We present the main result of this paper. First, we start with a factor $P$ of type ${\rm III_1}$. Its {\it core} $P \rtimes_{\s^\omega} \R$ is a factor of type ${\rm II_\infty}$ and does not depend on the choice of the state (or weight) $\omega$. Consider also a free, ergodic, measure-preserving flow $\gamma : \R \curvearrowright (X, \mu)$. The object of study is essentially the von Neumann algebra $(L^{\infty}(X) \otimes P) \rtimes_{\gamma \otimes \s^\omega} \R$. This von Neumann algebra turns out to be a factor of type ${\rm III_1}$ whose core is canonically isomorphic to $(L^{\infty}(X) \rtimes_{\g} \R) \otimes (P \rtimes_{\s^\omega} \R)$ (see Theorem $\ref{construction}$). Under this identification, the {\it dual action} is given by $\theta_s = \widehat{\gamma}_{-s} \otimes \theta^\omega_s$, with $\widehat{\gamma}$ the dual action of $\gamma$ and $\theta^\omega$ the dual action of $\s^\omega$. Notice that since $\gamma$ is free, ergodic and measure-preserving, and $\R$ is amenable, $L^{\infty}(X) \rtimes_{\g} \R$ is the hyperfinite type ${\rm II_\infty}$ factor.

We shall prove a result of complete classification for this family of factors in the case the core $P \rtimes_{\s^\omega} \R$ is a {\it full} factor. For this, we shall use a very nice result by Popa on unique tensor product decomposition of McDuff ${\rm II_1}$ factors. The following theorem is at the heart of this paper:
\begin{theo}\cite{popaduff}
Let $R_1$ and $R_2$ be two copies of the hyperfinite type ${\rm II_1}$ factor; let $N_1$ and $N_2$ be two full factors of type ${\rm II_1}$. Let us assume that the factor $N$ has both of the following decompositions: $N = R_1 \otimes N_1 = R_2 \otimes N_2$. Then there exist $t > 0$ and a unitary $u \in N$ such that $R_2 = u^*(R_1)^{1/t}u$ and $N_2 = u^*(N_1)^tu$. 
\end{theo}
We can give an equivalent version of this theorem in the type ${\rm II_\infty}$ case:
\begin{theo}[Mc Duff ${\rm II_\infty}$ factors]\label{duffinfini}
Let $R_1^{\infty}$ and $R_2^{\infty}$ be two copies of the hyperfinite type ${\rm II_{\infty}}$ factor; let $N_1^{\infty}$ and $N_2^{\infty}$ be two full factors of type ${\rm II_{\infty}}$. Let $\a : R_1^{\infty} \otimes N_1^{\infty} \to R_2^{\infty} \otimes N_2^{\infty}$ be an isomorphism. Then, there exist a unitary $u \in R_2^{\infty} \otimes N_2^{\infty}$ and two isomorphisms $\b : R_1^{\infty} \to R_2^{\infty}$, $\g : N_1^{\infty} \to N_2^{\infty}$  such that for any $x \in R_1^{\infty} \otimes N_1^{\infty}$, $u^*\a(x)u = (\b \otimes \g)(x)$.
\end{theo}



We can now give our main result concerning the complete classification of the construction in the case $P\rtimes_{\s^\omega} \R$ is a full factor.
\begin{theo}\label{mainresult}
Let $P$ be the factor of type ${\rm III_1}$ whose core is isomorphic to $L(\mathbf{F}_{\infty}) \otimes B(\h)$ and the dual action is given by the trace-scaling one-parameter automorphism group $(\a_t)_{t \in \R}$ from R\u{a}dulescu \cite{radulescu1991}. Let $(X_i, \mu_i, \g^i)$ $(i = 1, 2)$ be two measure spaces endowed with measure-preserving, free, ergodic flows. Denote $M_i = (L^{\infty}(X_i) \otimes P) \rtimes_{\gamma^i \otimes \s^\omega} \R$ $(i = 1, 2)$. Then $M_1$ and $M_2$ are isomorphic if and only if $\g^1$ and $\g^2$ are conjugate.
\end{theo}

Since the construction is canonical, it is obvious that if $\gamma^1$ and $\gamma^2$ are conjugate then $M_1$ and $M_2$ are isomorphic. Let us give a few ideas about the proof of the converse. Assume that $M_1$ and $M_2$ are isomorphic. Then, we know that the cores of $M_1$ and $M_2$ are isomorphic and the dual actions are cocycle conjugate (see \cite{{takesakiII}, {takesaki73}} for further details). But for $i = 1, 2$, the core of $M_i$, denoted by $\Core(M_i)$, is isomorphic to the McDuff type ${\rm II_\infty}$ factor $(L^{\infty}(X_i) \rtimes_{\g^i} \R) \otimes (P \rtimes_{\s^\omega} \R)$. If we apply the result by Popa (Theorem $\ref{duffinfini}$), we get on the \textquotedblleft first leg\textquotedblright of the tensor product decomposition of $\Core(M_1)$ and $\Core(M_2)$ an isomorphism $\pi : L^{\infty}(X_1) \rtimes_{\g^1} \R \to L^{\infty}(X_2) \rtimes_{\g^2} \R$ and a family $(v_s)_{s \in \R}$ of unitaries in $\Core(M_1)$ such that for any $z \in L^{\infty}(X_1) \rtimes_{\g^1} \R$
$$\pi^{-1}\left(\widehat{\g^2}_{s}(\pi(z))\right) \otimes 1 = v_s(\widehat{\g^1}_{s}(z) \otimes 1)v_s^*.$$
Using classical techniques, we show that there exists a family of unitaries $(w_s)_{s \in \R}$ which is a one-cocycle for $\gamma_1$ and such that
$$\pi^{-1}\widehat{\g^2}_{s}\pi = w_s\widehat{\g^1}_{s}w_s^*.$$
At last, using classical results by Takesaki \cite{{takesakiII}, {takesaki73}}, we prove that $\gamma^1$ and $\gamma^2$ are necessarily conjugate. Consequently, if the factor $P$ satisfies the assumptions of Theorem $\ref{mainresult}$ (and more generally if the core of $P$ is full), the factor of type ${\rm III_1}$, $(L^{\infty}(X) \otimes P) \rtimes_{\gamma \otimes \s^\omega} \R$, entirely remembers the flow $\gamma$.

Finally, when $P$ is a free Araki-Woods factor \cite{shlya97}, we will compute Connes' $\tau$ invariant for the factor $(L^{\infty}(X) \otimes P) \rtimes_{\gamma \otimes \s^\omega} \R$. We will exhibit a new example of an uncountable family of pairwise non-isomorphic factors of type ${\rm III_1}$ with the same $\tau$ invariant.

\section{The General Construction}

\paragraph{\bf Notation}
We want to introduce here all the notations we will use in this paper. Let $(A, (\a_t)_{t \in \R})$ and $(B, (\b_t)_{t \in \R})$ be two von Neumann algebras endowed with a one-parameter automorphism group. Let $\mathcal{F} : L^2(\R) \to L^2(\R)$ be the Fourier transform on $\R$: for $f \in L^1(\R)\cap L^2(\R)$, $\mathcal{F}(f) = \xi \mapsto \int_{\R}f(t)\exp(-it\xi)dt$. So, if $\rho = \Ad(\mathcal{F}) : L^{\infty}(\R) \to B(L^2(\R))$, we get for any $t \in \R$,  $\rho(e^{it \cdot}) = \rho_t$, with $\rho_t$ the translation operators on $L^2(\R)$.  Let $\a : A \to A \otimes L^{\infty}(\R)$ defined by: for any $x \in N$ and any $t \in \R$, $\a(x)(t) = \a_{-t}(x)$. Let $\b : B \to L^{\infty}(\R) \otimes B$ defined by: for any $y \in B$ and any $t \in \R$, $\b(y)(t) = \b_{-t}(y)$. 
\begin{df}\label{defi}
We define the \emph{crossed product} of $A$ and $B$ under $\a$ and $\b$ in the following way: it is the von Neumann subalgebra of $A \otimes B(L^2(\R)) \otimes B$ generated by $(\a(A) \otimes 1 \cup 1 \otimes \tilde{\b}(B))$, with $\tilde{\b} = (\rho \otimes 1)\b$. We shall denote it by $A \indice{\Join}_{\b} B$.
\end{df} 
It is clear that $A \indice{\Join}_{\b} B$ is canonically isomorphic to the von Neumann subalgebra of $A \otimes B(L^2(\R)) \otimes B$ generated by $(\tilde{\a}(A) \otimes 1 \cup 1 \otimes \b(B))$, with $\tilde{\a} = (1 \otimes \rho)\a$. Let $(\gamma_t)_{t \in \R}$ be a free, ergodic flow on $(X, \mu)$ and let $(N, (\alpha_t)_{t \in \R})$ be a factor of type ${\rm II_{\infty}}$ endowed with a trace-scaling one-parameter automorphism group. Let $P$ be the factor of type ${\rm III_1}$, $N \rtimes_{\a} \R$. Let $\tau_N$ be a faithful, semi-finite, normal trace on $N$ (this trace is unique up to a scalar $\la > 0$) and let $\omega$ be the dual weight on $P$ of $\tau_N$ under the action $\a$; let $\s^{\omega}$ be the modular automorphism group of the weight $\omega$ and $\b$ the action $\g \otimes \s^{\omega}$ of $\R$ on $L^{\infty}(X)\otimes P$. We shall denote by $C$ the crossed product $N \indice{\Join}_{\g} L^{\infty}(X)$ of $N$ and $L^{\infty}(X)$ under $\a$ and $\g$. We will prove in Proposition $\ref{identification}$ that $C \otimes B(L^2(\R))$ is canonically isomorphic to $(L^{\infty}(X) \otimes P) \rtimes_{\g \otimes \s^{\omega}} \R$. We shall denote by $\phi$ the n.f.s. weight $\tau \otimes \omega$ on $L^{\infty}(X)\otimes P$, and $\tilde{\phi}$ the dual weight of $\phi$ on $(L^{\infty}(X) \otimes P) \rtimes_{\g \otimes \s^{\omega}} \R$ under the action $\beta$. Moreover, most of the time $L^{\infty}(X)\otimes P \otimes B(L^2(\R))$ will be denoted by $M$; the canonical embedding of $(L^{\infty}(X) \otimes P) \rtimes_{\g \otimes \s^{\omega}} \R$ into $M$ will be denoted by $\pi$. At last we know, according to Proposition $8.4$ in \cite{takesaki73}, that $(N \rtimes_{\a} \R) \rtimes_{\s^{\omega}} \R$ is canonically isomorphic to $N \otimes B(L^2(\R))$, i.e. $P \rtimes_{\s^{\omega}} \R = N \otimes B(L^2(\R))$. We shall always identify both of these factors.

The following proposition gives a justification to Definition $\ref{defi}$.
\begin{prop}\label{identification}
Let $(A, (\a_t)_{t \in \R})$ and $(B, (\b_t)_{t \in \R})$ be two von Neumann algebras endowed with a one-parameter automorphism group. Let $C = A \indice{\Join}_{\b} B$.
\begin{enumerate}
\item If $B = L^{\infty}\left(\widehat{\Gamma}\right)$, with $\Gamma$ subgroup of $\R$ and $(\b_t)_{t \in \R}$ is the action of $\R$ by translations, then $C$ is canonically isomorphic to $A \rtimes_{\a} \G$.

\item If $A = A_0 \rtimes_{\a^0} \R$ and $\a$ is the dual action of $\a^0$, then $C$ is canonically isomorphic to $(A_0 \otimes B) \rtimes_{\a^0 \otimes \b} \R$.
\end{enumerate}
\end{prop}
\begin{proof}
$(1)$ Let $B = L^{\infty}\left(\widehat{\Gamma}\right)$, with $(\b_t)_{t \in \R}$ the action of $\R$ by translations. Let $G = \widehat{\G}$. For $\g \in \G$, let $u_{\g} \in L^{\infty}(G)$ defined by $u_{\g}(g) = \left\langle g, \g\right\rangle$. Then for any $t \in \R$,
\begin{eqnarray*}
\b_t(u_{\g})(g) & = & \langle g - t, \g\rangle \\
& = & \exp(-it \g)\langle g, \g\rangle \\
& = & \exp(-it \g)u_{\g}(g).
\end{eqnarray*}
Thus, for any $\gamma \in \Gamma$, $u_{\gamma}$ is an eigenvector associated with $\b_t$ for the eigenvalue $\exp(-it\gamma)$. Therefore, $\tilde{\b}(u_{\g}) = \rho_{-\g} \otimes u_{\g}$ and consequently
\begin{eqnarray*}
(1 \otimes \tilde{\b}(u_{\gamma}))(\a(a) \otimes 1)(1 \otimes \tilde{\b}(u_{\gamma}))^* & = & (1 \otimes \rho_{-\gamma})\a(a)(1 \otimes \rho_{-\gamma})^* \otimes 1\\
& = & \a(\a_{-\gamma}(a)) \otimes 1.
\end{eqnarray*}
As $L^{\infty}(G)$ is spanned by the operators $u_{\g}$ for $\g \in \G$, it is clear that $A \indice{\Join}_{\b} B$ is canonically isomorphic to $A \rtimes_{\a} \G$.

$(2)$ Let $A = A_0 \rtimes_{\a^0} \R$ and let $\a$ be the dual action of $\a^0$. We still denote by $\a^0 : A_0 \to A$, the mapping defined by: for any $x \in A_0$ and any $t \in \R$, $\a^0(x)(t) = \a^0_{-t}(x)$. We shall denote by $\la_t \in A$ the unitaries which implement the action $\a^0$ of $\R$ on $A_0$. By definition of the dual actions \cite{takesaki73}, for any $a \in A_0$ and any $t \in \R$, we get: 
\begin{eqnarray*}
\tilde{\a}(\a^0(a)) & = & \a^0(a) \otimes 1 \\
\tilde{\a}(\la_t) & = & \la_t \otimes \rho_t.
\end{eqnarray*}
The von Neumann algebra $A \indice{\Join}_{\b} B$ is generated by $\tilde{\a}(\a^0(A_0))\otimes1$, $\tilde{\a}(\la_t) \otimes1$ for $t \in \R$ and $1 \otimes \b(B)$. But $\tilde{\a}(A_0)\otimes1$ is isomorphic to $A_0$, $1 \otimes \b(B)$ is isomorphic to $B$, and $\tilde{\a}(A_0)\otimes1$ and $1 \otimes \b(B)$ commute in $A \indice{\Join}_{\b} B$. Moreover, for any $a \in A_0$ and any $b \in B$, we have
\begin{eqnarray*}
\tilde{\a}(\la_t)\tilde{\a}(\a^0(a))\tilde{\a}(\la_t)^* & = & \a^0(\a^0_t(a)) \\
(\tilde{\a}(\la_t)\otimes1)(1\otimes\b(b))(\tilde{\a}(\la_t)\otimes1)^* & = & 1\otimes\b(\b_t(b)).
\end{eqnarray*}
Finally, after trivial identifications, we get that $A \indice{\Join}_{\b} B$ is canonically isomorphic to $(A_0 \otimes B) \rtimes_{\a^0 \otimes \b} \R$.
\end{proof}

From now on, let $(N, \a)$ be a factor of type ${\rm II_{\infty}}$ endowed with a trace-scaling one-parameter automorphism group and let $(\g_t)_{t \in \R}$ be a finite or infinite measure-preserving, free, ergodic flow on the measure space $(X, \mu)$. Let $P = N \rtimes_{\a} \R$ and $\b = \g \otimes \s^{\omega}$. According to the theorem of duality by Takesaki, we know that $P \rtimes_{\s^{\omega}} \R$ is canonically isomorphic to $N \otimes B(L^2(\R))$. Thus, we get thanks to Proposition $\ref{identification}$ that $C \otimes B(L^2(\R))$ is canonically isomorphic to $(L^{\infty}(X) \otimes P) \rtimes_{\b} \R$, with $C$ the crossed product $N \indice{\Join}_{\g} L^{\infty}(X)$. From now on, we shall identify $C \otimes B(L^2(\R))$ with $(L^{\infty}(X) \otimes P) \rtimes_{\b} \R$. For any weight $\psi$, we shall denote by $\Delta_{\psi}$ the modular operator associated with $\psi$. Let $\tau$ be the faithful, normal, semifinite trace on $L^{\infty}(X)$ given by the measure $\mu$ (preserved by the flow $(\g_t)_{t \in \R}$). It is clear that $\phi = \tau \otimes \omega$ is a faithful, normal, semifinite weight on $L^{\infty}(X) \otimes P$. As $\tau$ is preserved by the action $\g_t$ and $\omega$ by $\s^{\omega}_t$, it is obvious that $\phi$ is preserved by the action $\b_t = \g_t \otimes \s^{\omega}_t$. Let $\tilde{\phi}$ the dual weight of $\phi$ on $(L^{\infty}(X) \otimes P) \rtimes_{\b} \R$ under the action $\b$ and $\la_s$ the unitaries which implement the action $\b$ of $\R$ on $L^{\infty}(X) \otimes P$. According to Proposition $5.15$ in \cite{takesaki73}, we know that $\s^{\tilde{\phi}}$ acts on $(L^{\infty}(X) \otimes P) \rtimes_{\b} \R$ in the following way: for any $x \in L^{\infty}(X) \otimes P$ and $s \in \R$,
\begin{eqnarray*}
\s_t^{\tilde{\phi}}(\pi_{\b}(x)) & =  & \pi_{\b}(\s_t^{\phi}(x))\\
\s_t^{\tilde{\phi}}(\la_s) & = & \la_s.
\end{eqnarray*}
Let us denote by $L^2(X, \mu)$ and $L^2(P, \omega)$ the Hilbert spaces on which  $L^{\infty}(X)$ and $P$ act canonically. Hence, the Hilbert space on which $(L^{\infty}(X) \otimes P) \rtimes_{\b} \R$ acts is nothing but $L^2(X) \otimes L^2(P, \omega) \otimes L^2(\R)$; by definition of the dual weight $\tilde{\phi}$ on $(L^{\infty}(X) \otimes P) \rtimes_{\b} \R$ (see \cite{takesaki73} for further details), we have 
\begin{eqnarray*}
\forall t \in \R, \: \Delta_{\tilde{\phi}}^{it} & = & \Delta_{\phi}^{it} \otimes 1\\
                                                                    & = & 1 \otimes \Delta_{\omega}^{it} \otimes 1.
\end{eqnarray*}
On the other hand, as the family $(\la_t)$ is in the centralizer of the weight $\tilde{\phi}$, it is a one-cocycle for the action $\s^{\tilde{\phi}}$. We know according to Theorem $1.2.4$ of \cite{connes73} that there exists a faithful, normal, semifinite weight $\psi$ on $(L^{\infty}(X) \otimes P) \rtimes_{\b} \R$ such that:
$$\forall t \in \R, \: \s_t^{\psi} = \la_t^*\s_t^{\tilde{\phi}}\la_t.$$
We denote by $U$ the unitary representation of $\R$ on $L^2(X, \mu)$ which implement the action $\g$ of $\R$ on $L^{\infty}(X)$. For any $x \in L^{\infty}(X) \otimes P$, we have
\begin{eqnarray*}
\s_t^{\psi}(\pi_{\b}(x)) & = & \la_t^*\s_t^{\tilde{\phi}}(\pi_{\b}(x))\la_t\\
 & = &  \la_t^*\pi_{\b}(\s_t^{\phi}(x))\la_t\\
 & = & \pi_{\b}(\b_{-t}(\s_t^{\phi}(x)))\\
 & = & \pi_{\b}((\g_{-t} \otimes \id)(x)).
\end{eqnarray*}
Hence, by definition of the modular operator,
\begin{equation}\label{presqueper}
\forall t \in \R, \: \Delta_{\psi}^{it} = U_{-t} \otimes 1 \otimes 1.
\end{equation}
Thus, from the Equation $(\ref{presqueper})$, it is not difficult to see that, since the action $\g$ is ergodic, the centralizer of $\psi$ (denoted by $\Cent(\psi)$) is the von Neumann algebra $(\pi_{\s^{\omega}}(P)\cup\la(\R))''$ spanned by $\pi_{\s^{\omega}}(P)$ and $\la(\R)$; after trivial identifications, it is nothing but $N \otimes B(L^2(\R))$ and consequently $\Cent(\psi)$ is a factor of type ${\rm II_{\infty}}$. We are now able to prove the following theorem:
\begin{theo}
Let $(N, \a)$ be a factor of type ${\rm II_{\infty}}$ endowed with a trace-scaling one-parameter automorphism group and $(\g_t)_{t \in \R}$ be a finite or infinite measure-preserving, free, ergodic flow on $(X, \mu)$. Then the crossed product $C = N \indice{\Join}_{\g} L^{\infty}(X)$ is a factor of type ${\rm III_1}$.
\end{theo}
\begin{proof}
Let us prove first that $C$ is a factor. We remind that $C \otimes B(L^2(\R))$ is canonically isomorphic to $(L^{\infty}(X) \otimes P) \rtimes_{\b} \R$. We denote by $\pi$ the natural embedding of $(L^{\infty}(X) \otimes P) \rtimes_{\b} \R$ into $M =  L^{\infty}(X) \otimes P \otimes B(L^2(\R))$. As $P$ is a factor of type ${\rm III_1}$, we know according to Theorem ${\rm XII}$ $1.7$ of \cite{takesakiII} that
\begin{equation*}
\{P \otimes B(L^2(\R))\} \cap (P \rtimes_{\s^{\omega}} \R)' = \C.
\end{equation*}
So, it is not very difficult to see that
\begin{equation*}
\{(L^{\infty}(X) \otimes P) \rtimes_{\g \otimes \s^{\omega}} \R\} \cap (1 \otimes \pi(P \rtimes_{\s^{\omega}} \R))' \subset L^{\infty}(X)^{\g} \otimes1\otimes1,
\end{equation*}
with $L^{\infty}(X)^{\g}$ the fixed points subalgebra of $L^{\infty}(X)$ under the action $\g$. But, as $\g$ is ergodic, we know that $L^{\infty}(X)^{\g} = \C$. Consequently, we obtain that
\begin{equation}\label{ergodic}
\{(L^{\infty}(X) \otimes P) \rtimes_{\g \otimes \s^{\omega}} \R\} \cap (1 \otimes \pi(P \rtimes_{\s^{\omega}} \R))' = \C
\end{equation}
and $C$ is a factor.

Furthermore, we want to prove that $C$ is a factor of type ${\rm III}$. Indeed, suppose that $C$ were semifinite. Then $(L^{\infty}(X) \otimes P) \rtimes_{\b} \R$ would be semifinite and the modular automorphism group $\s^{\psi}_t$ would be inner. Hence, there would exist unitaries $u_t \in (L^{\infty}(X) \otimes P) \rtimes_{\b} \R$ such that $\s^{\psi}_t = \Ad(u_t)$ for any $t \in \R$. But, this implies that $u_t \in (L^{\infty}(X) \otimes P) \rtimes_{\b} \R\cap\Cent(\psi)'$, and thanks to Equation $(\ref{ergodic})$, we get $\s^{\psi}_t = \Id$ for any $t \in \R$. That means exactly that $\psi$ would be a trace and consequently $\Delta_{\psi} = \Id$. Thus, according to Equation $(\ref{presqueper})$, $U_t = 1$ for any $t \in \R$ which contradicts the fact that $(\g_t)$ is free. Therefore, $C$ is a factor of type ${\rm III}$. Moreover, as $\Cent(\psi)$ is a factor, we know according to Corollary $3.2.7$ of \cite{connes73}, that $C$ cannot be a factor of type ${\rm III_0}$. So now, if we want to find out the under-type of $C$, we can compute $T$ invariant of Connes \cite{connes73}. We remind that we have showed that $\s_t^{\psi}(\pi_{\b}(x)) = \pi_{\b}((\g_{-t} \otimes \id)(x))$ for any $x \in L^{\infty}(X) \otimes P$. As the flow $(\g_t)$ is supposed to be free, it cannot be periodic; moreover, for the same reasons as before, $\s^{\psi}_t$ cannot be inner. Therefore $T(C) = \{0\}$ and $C$ is a factor of type ${\rm III_1}$.
\end{proof}

Before proceeding, we are going to remind the definition of a strongly ergodic action of a locally compact group, the definition of a full factor and state a well-known theorem.
\begin{df}\cite{{schmidt1980}, {schmidt1981}}
Let $G$ be a locally compact group which acts on a probability space $X$ by $(\g_g)$ in a finite measure preserving way. This action is said to be \emph{strongly ergodic} when for any sequence of projections $(p_n)$ in $L^{\infty}(X)$, if the sequence $(\g_g(p_n) - p_n)$ tends to $0$ $\ast-$strongly uniformly on compacts sets, then the sequence $(p_n - \tau(p_n)1)$ tends to $0$ $\ast-$strongly.
\end{df}
\begin{df}\cite{connes74}
Let $C$ be a factor. Let $(x_n)_{n \in \N}$ be a bounded sequence in $C$. The sequence $(x_n)$ is said to be \emph{centralising} if for any normal state $\phi$ on $C$, $\|[x_n, \phi]\| \to 0$ when $n \to +\infty$; it is said to be \emph{trivial} if there exists a sequence of complex numbers $(\la_n)_{n \in \N}$ such that $x_n - \la_n \to 0$ $*-$strongly. At last, the factor $C$ is said to be a \emph{full factor} if any centralising sequence in $C$ is trivial.
\end{df}
\begin{theo}\cite{{schmidt1980}, {schmidt1981}}
Let $G$ be a locally compact group. If $G$ has Kazhdan property $(T)$, then any ergodic action of $G$ on a probability space is strongly ergodic. If $G$ is amenable, then it cannot act strongly ergodically on a probability space.
\end{theo}

We are now able to state the following general proposition: 
\begin{prop}\label{propgen}
Let $(N, \a)$ be a factor of type ${\rm II_{\infty}}$ endowed with a trace-scaling one-parameter automorphism group and $(\g_t)_{t \in \R}$ be a finite measure-preserving, free, ergodic flow on $(X, \mu)$. Then the crossed product $C = N \indice{\Join}_{\g} L^{\infty}(X)$ is a non-full factor.
\end{prop}
\begin{proof}
We still denote by $M$ the von Neumann algebra $L^{\infty}(X) \otimes P \otimes B(L^2(\R))$ and we remind that $C \otimes B(L^2(\R))$ is canonically isomorphic to the factor $(L^{\infty}(X) \otimes P) \rtimes_{\g \otimes \s^{\omega}} \R$. We remind that $\pi$ is the embedding of $(L^{\infty}(X) \otimes P) \rtimes_{\b} \R$ into $M$. In order to prove that $C$ is not a full factor, we have thus to find a centralising sequence which is not trivial. The action $(\g_t)$ of $\R$ on $L^{\infty}(X)$ is not strongly ergodic because $\R$ is abelian and thus amenable. Therefore, there exists a sequence of projections $(p_n)$ in $L^{\infty}(X)$ such that $(\g_t(p_n) - p_n)$ tends to $0$ $\ast-$strongly, uniformly on compact sets, but $(p_n - \tau(p_n)1)$ does not tend to $0$ $\ast-$strongly. We have now to study the behaviour of $(p_n)$ in $(L^{\infty}(X) \otimes P) \rtimes_{\b} \R$. The fact that $(\g_t(p_n) - p_n)$ tends to $0$ $\ast-$strongly, uniformly on compact sets implies that in $M$, $(\pi(p_n) - p_n\otimes1\otimes 1)$ tends to 0 $\ast-$strongly. Furthermore, every single element $p_n\otimes1 \otimes 1$ is central in $M$ and therefore it commutes trivially with any normal form on $M$: hence $(p_n\otimes1 \otimes 1)$ is a centralising sequence in $M$. As $(\pi(p_n) - p_n\otimes1 \otimes 1)$ tends to 0 $\ast-$strongly, $(\pi(p_n))$ is also a centralising sequence in $M$. Consequently, $(p_n)$ is also a centralising sequence in $(L^{\infty}(X) \otimes P) \rtimes_{\b} \R$. Moreover, one can easily prove that $(p_n)$ is not a trivial sequence in $(L^{\infty}(X) \otimes P) \rtimes_{\b} \R$ because $(p_n - \tau(p_n)1)$ does not tend to $0$ $\ast-$strongly. Consequently, $C$ is not full.
\end{proof}
We can notice that $C$ may be full when the preserved measure $\mu$ is infinite. For example, if $\R$ acts by translations on $(\R, \la)$, then $N \indice{\Join}_{\g} L^{\infty}(\R)$ is isomorphic to $P = N \rtimes_{\a} \R$ according to Proposition \ref{identification}, which is full if $N$ is full \cite{shlya2004}. We are now going to give a precise description of the core of the factor $C$.
\begin{theo}\label{construction}
Let $(N, \a)$ be a factor of type ${\rm II_{\infty}}$ endowed with a trace-scaling one-parameter automorphism group and $(\g_t)_{t \in \R}$ be a finite or infinite measure-preserving, free, ergodic flow on $(X, \mu)$. Then the crossed product $C = N \indice{\Join}_{\g} L^{\infty}(X)$ is a factor of type ${\rm III_1}$ whose core is isomorphic to
$$(L^{\infty}(X) \rtimes_{\g} \R) \otimes N.$$
Under this identification, the dual action of $\R$ on the core of $C$, denoted by $(\t_s)$, is given by $\t_s = \widehat{\g}_{-s} \otimes \a_s$, with $\widehat{\g}$ the dual action of $\g$. 
\end{theo} 
\begin{proof}
For the convenience of the proof, we shall denote by $N_1$ the crossed product $P \rtimes_{\s^{\omega}} \R$; we know that $N_1$ is nothing but $N \otimes B(L^2(\R))$, and according to Theorem $4.6$ in \cite{takesaki73}, the dual action of $\s^{\omega}$ on $N_1$ is nothing but $\a_t \otimes \Ad(\rho^*_t)$ with $\rho_t$ the left regular representation of $\R$ on $B(L^2(\R))$. We remind that $\tau$ is the faithful, normal, semifinite trace on $L^{\infty}(X)$ given by the measure $\mu$ (preserved by the flow $(\g_t)_{t \in \R}$). It is clear that $\phi = \tau \otimes \omega$ is a faithful, normal, semifinite weight on $L^{\infty}(X) \otimes P$. As $\tau$ is preserved by the action $\g_t$ and $\omega$ by $\s^{\omega}_t$, it is obvious that $\phi$ is preserved by the action $\b_t = \g_t \otimes \s^{\omega}_t$. Let $\tilde{\phi}$ be the dual weight of $\phi$ on $(L^{\infty}(X) \otimes P) \rtimes_{\b} \R$ under the action $\b$. Let us denote by $\h = L^2(X, \mu)$ et $\k = L^2(P, \omega)$ the Hilbert spaces on which $L^{\infty}(X)$ and $P$ act canonically. Hence, $(L^{\infty}(X) \otimes P) \rtimes_{\b} \R$ acts on $L^2(\R, \h \otimes \k)$. Let us denote by $\Core((L^{\infty}(X) \otimes P) \rtimes_{\b} \R)$ the core of the factor $(L^{\infty}(X) \otimes P) \rtimes_{\b} \R$. By definition, we have $\Core((L^{\infty}(X) \otimes P) \rtimes_{\b} \R) = ((L^{\infty}(X) \otimes P) \rtimes_{\b} \R) \rtimes_{\s^{\tilde{\phi}}} \R$. The core of $(L^{\infty}(X) \otimes P) \rtimes_{\b} \R$ acts on the Hilbert space $L^2(\R \times \R, \h \otimes \k)$, and it is spanned by $\pi_{\s^{\tilde{\phi}}}(N)$ and $\rho(\R)$, if we denote by $\rho_s$ the unitaries which implement the action $\s^{\tilde{\phi}}$ of $\R$ on $(L^{\infty}(X) \otimes P) \rtimes_{\b} \R$. Therefore, $\Core((L^{\infty}(X) \otimes P) \rtimes_{\b} \R)$ is spanned by three types of operators:
\begin{enumerate}
\item $\widetilde{a \otimes b} = \pi_{\s^{\tilde{\phi}}} \circ \pi_{\b}(a \otimes b)$, for $a \in L^{\infty}(X)$ and $b \in P$; for any $\zeta \in L^2(\R \times \R, \h \otimes \k)$, we have $((\widetilde{a \otimes b})\zeta)(s', t') = \left(\g_{-t'}(a) \otimes \s_{-(s' + t')}^{\omega}(b)\right)\zeta(s', t')$.
\item $u_t = \pi_{\s^{\tilde{\phi}}}(\la_t)$, for $t \in \R$; for any $\zeta \in L^2(\R \times \R, \h \otimes \k)$, we have $(u_t\zeta)(s', t') = \zeta(s', t' - t)$.
\item $v_s = \rho_s$, for $s \in \R$; for any $\zeta \in L^2(\R \times \R, \h \otimes \k)$, we have $(v_s\zeta)(s', t') = \zeta(s' - s, t')$.
\end{enumerate}
Let $U$ be the following unitary: 
$$U : \fonction{L^2(\R, \h) \otimes L^2(\R, \k)}{L^2(\R \times \R, \h \otimes \k)}{\xi \otimes \eta}{\left\{(s', t') \mapsto \xi(t') \otimes \eta(s' + t')\right\}}.$$
Let us denote by $\pi_{\g}(L^{\infty}(X))$ and $\la^g(\R)$ the images of $L^{\infty}(X)$ and $\R$ in the crossed product $L^{\infty}(X, \mu) \rtimes_{\g} \R$; in the same way, we denote by $\pi_{\s^{\omega}}(P)$ and $\la^d(\R)$ the images of $P$ and $\R$ in $P \rtimes_{\s^{\omega}} \R = N_1$. But now the question is: what is the behaviour of our three types of operators when we intertwine them by the unitary $U$? Let $\xi \in L^2(\R, \h)$, $\eta \in L^2(\R, \k)$, $s$, $s'$, $t$ and $t' \in \R$. We have
\begin{eqnarray*}
\left(U^*\widetilde{a \otimes b}U(\xi \otimes \eta)\right)(s', t') & = & \g_{-t'}(a)\xi(t') \otimes \s_{-s'}^{\omega}(b)\eta(s')\\
\left(U^*u_tU(\xi \otimes \eta)\right)(s', t') & = & \xi(t' - t) \otimes \eta(s' - t)\\
\left(U^*v_sU(\xi \otimes \eta)\right)(s', t') & = & \xi(t') \otimes \eta(s' - s).
\end{eqnarray*}
Using notations we have just introduced, we easily obtain: 
\begin{eqnarray}\label{un}
U^*\widetilde{a \otimes b}U & = & \pi_{\g}(a) \otimes \pi_{\s^{\omega}}(b)\\ \label{deux}
U^*u_tU & = & \la_t^g \otimes \la_t^d\\ 
U^*v_sU & = & 1 \otimes \la_s^d. \label{trois}
\end{eqnarray}
Therefore at this stage, we have proved that the core of $(L^{\infty}(X) \otimes P) \rtimes_{\b} \R$, intertwined by $U$, is equal to
$$\left(L^{\infty}(X) \rtimes_{\g} \R\right) \otimes N_1.$$
It remains us to understand what the dual action of $\s^{\tilde{\phi}}$ (denoted by $\t_s$) is. According to  \cite{takesaki73}, for any $g \in \R$ and $s$, $t \in \R$, we have:
\begin{eqnarray*}
\theta_g(\widetilde{a \otimes b}) & = & \widetilde{a \otimes b}\\ 
\theta_g(u_t) & = & u_t\\ 
\theta_g(v_s) & = & e^{-isg}v_s. 
\end{eqnarray*}
For any $g \in \R$, let $\theta'_g = U^*\theta_gU$;  we notice, thanks to Equations $(\ref{un})$, $(\ref{deux})$ and $(\ref{trois})$ that: 
\begin{eqnarray*}
\theta'_g\left(\pi_{\g}(a) \otimes \pi_{\s^{\omega_0}}(b)\right) & = & \pi_{\g}(a) \otimes \pi_{\s^{\omega_0}}(b)\\
\theta'_g\left(\la_t^g \otimes 1\right) & = & e^{itg}\left(\la_t^g \otimes 1\right)\\
\theta'_g\left(1 \otimes \la_s^d\right)& = & e^{-isg}\left(1 \otimes \la_s^d\right).
\end{eqnarray*}
Finally, as $\a_t \otimes \Ad(\rho^*_t)$ is the dual action of $\s^{\omega}$, we know according to \cite{takesaki73} that the dual action $\t'_g$ is given by:
$$\forall g \in \R, \; \theta'_g = \widehat{\g}_{-g} \otimes \a_g \otimes \Ad(\rho^*_g).$$
As $C \otimes B(L^2(\R))$ is canonically isomorphic to $(L^{\infty}(X) \otimes P) \rtimes_{\b} \R$, after trivial identifications, we obtain that the core of $C$ is nothing but $(L^{\infty}(X) \rtimes_{\g} \R) \otimes N$ and the dual action is given by $\widehat{\g}_{-s} \otimes \a_s$.
\end{proof}

\section{Classification in the Case $N$ is a Full Factor} 

In this section, we shall assume that $N$ is a full factor of type ${\rm II_{\infty}}$ endowed with a trace-scaling one-parameter automorphism group denoted by $(\a_t)_{t \in \R}$ and $(\g_t)_{t \in \R}$ is a free, ergodic, measure-preserving flow on $(X, \mu)$. The first example of such a full factor $N$ has been given by R\u{a}dulescu in \cite{radulescu1991}: it is $L(\F_{\infty}) \otimes B(\h)$. We can notice that other examples exist thanks to the theory of free Araki-Woods factors developed by Shlyakhtenko in \cite{{shlya2003}, {shlya2004}, {shlya97}}.
 
Before stating the main theorem of this section, let us introduce a few notations. Let $(X, \mu)$ and $(X', \mu')$ be two measure spaces on which $\R$ acts by $\g$ and $\g'$ in a free, ergodic and measure-preserving way. Let $(N, \a)$ and $(N', \a')$ be two full factors of type ${\rm II_{\infty}}$, both of them endowed with a trace-scaling one-parameter automorphism group denoted by $\a$ and $\a'$. According to \cite{takesaki73}, the covariant systems $(N, \a)$ and $(N', \a')$ are said to be \emph{weakly equivalent}, if there exist an isomorphism $\pi : N \to N'$ and a strongly continuous family of unitaries $(u_s)$ of $N$ such that $(u_s)$ is a one-cocycle for the action $\a$ and for any $x \in N$ and any $s \in \R$,
$$\pi^{-1}\a'_s\pi(x) = u_s\a_s(x)u_s^*.$$
Takesaki proved in \cite{takesaki73}, Corollary $8.4$, that it is equivalent to say that the factors $N \rtimes_{\a} \R$ and $N' \rtimes_{\a'} \R$ are isomorphic. 
Let us denote by $C$ the crossed product $N \indice{\Join}_{\g} L^{\infty}(X)$ and $C'$ the crossed product $N' \indicebis{\Join}_{\g'} L^{\infty}(X')$. As our construction is canonical in a obvious way, it is clear that if the flows $(X, \g)$ and $(X', \g')$ are conjugate and $(N, \a)$ and $(N', \a')$ are weakly equivalent, then $C$ and $C'$ are isomorphic (this claim is true even if $N$ and $N'$ are not full). We want to prove the converse here, when $N$ and $N'$ are full: 
\begin{theo}\label{bijection}
If $C$ and $C'$ are isomorphic, then $(N, \a)$ and $(N', \a')$ are weakly equivalent, and the flows $(X, \g)$ and $(X', \g')$ are conjugate, i.e. there exists an isomorphism $\pi : L^{\infty}(X) \to L^{\infty}(X')$ such that for any $t \in \R$, $\pi\g_t = \g'_t\pi$.
\end{theo}
First of all, we want to prove an important result about full factors of type ${\rm II_1}$.
\begin{lem}\label{inclusion}
Let $M$ be a full factor of type ${\rm II_1}$ and $N$ be any type ${\rm II_1}$ factor. Then, for any net of unitaries $(u_i)$ in $M \otimes N$ such that $\|xu_i - u_ix\|_2 \to 0$ for all $x \in M \otimes N$, 
\begin{equation*}
\sup_{x \in M, \|x\| \leq 1}\|u_i(x\otimes1) - (x\otimes1)u_i\|_2 \to 0.
\end{equation*}
\end{lem}
\begin{proof}
Let $J_M$ be the canonical anti-unitary on $L^2(M)$ associated with the trace $\tau_M$. Let $(u_i)$ be a net of unitaires in $M \otimes N$ such that $\|xu_i - u_ix\|_2 \to 0$ for all $x \in M \otimes N$. We are going to prove that $\|u_i - 1\otimes E_N(u_i)\|_2 \to 0$ with $E_N : M \otimes N \to N$ the conditional expectation defined by $E_N = \tau_M \otimes \Id$. It will prove our result; indeed, for any $x \in M$, as $x\otimes1$ commutes with $1\otimes E_N(u_i)$, we get $\|(x\otimes1)u_i - u_i(x\otimes1)\|_2 = \|(x\otimes1)v_i - v_i(x\otimes1)\|_2$, with $v_i = u_i - 1\otimes E_N(u_i)$.
Moreover, for any $x \in M$, $\|x\| \leq 1$, we have
\begin{eqnarray*}
\|(x\otimes1)v_i - v_i(x\otimes1)\|_2 & \leq & \|(x\otimes1)v_i\|_2 + \|v_i(x\otimes1)\|_2\\
& \leq & \|(x\otimes1)v_i\|_2 + \|(J_Mx^*J_M\otimes1)v_i\|_2\\
& \leq & \|x\| \: \|v_i\|_2 + \|J_Mx^*J_M\| \: \|v_i\|_2\\
& \leq & 2\|v_i\|_2.
\end{eqnarray*}
Consequently, we get
\begin{equation*}
\sup_{x \in M, \: \|x\| \leq 1}\|(x\otimes1)u_i - u_i(x\otimes1)\|_2 \leq 2\|u_i - 1\otimes E_N(u_i)\|_2 \to 0.
\end{equation*}
We shall denote by $\Omega_M \in L^2(M)$ and $\Omega_N \in L^2(N)$ the images of $1_M \in M$ and $1_N \in N$ associated with the GNS constructions respectively for $M$ and $N$. For any $x \in M$ and $x' \in M'$, $[u_i, xx' \otimes1] = (x'\otimes1)[u_i, x\otimes1]$. Therefore, as $(u_i)$ is a bounded net in $B(L^2(M)) \otimes N$, we get that for any $y \in C^*(M, M')$, $[u_i, y \otimes1] \to 0$ $*-$strongly. According to Theorem $2.1$ in \cite{connes76}, as $M$ is a full factor of type ${\rm II_1}$, we know that $K(L^2(M)) \subset C^*(M, M')$. Thus, if we denote by $P_{\Omega_M} \in B(L^2(M))$ the rank-one projection onto $\C\Omega_M$, we have proved that $[u_i, P_{\Omega_M} \otimes1] \to 0$ $*-$strongly.

We still denote by $\Omega_M : \C \to L^2(M)$ the map which sends $\la$ onto $\la\Omega_M$. It is easy to notice that $\Omega_M\Omega_M^{\ast} = P_{\Omega_M} : L^2(M) \to L^2(M)$.  We can notice that $E_N(u_i) = (\Omega_M^*\otimes1)u_i(\Omega_M\otimes1)$, thus \begin{eqnarray*}
E_N(u_i)^*E_N(u_i) & = & (\Omega_M^*\otimes1)u_i^*(P_{\Omega_M}\otimes1)u_i(\Omega_M\otimes1)\\
& = & (\Omega_M^*\otimes1)(P_{\Omega_M}\otimes1)(\Omega_M\otimes1) - (\Omega_M^*\otimes1)u_i^*[u_i, P_{\Omega_M}\otimes1](\Omega_M\otimes1).
\end{eqnarray*}
Finally, we get 
\begin{eqnarray*}
\|u_i - 1\otimes E_N(u_i)\|_2^2 & = &
\langle(u_i - 1\otimes E_N(u_i))(\Omega_M\otimes\Omega_N),(u_i - 1\otimes E_N(u_i))(\Omega_M\otimes\Omega_N)\rangle\\ 
& = & 1 - \tau_N(E_N(u_i)^*E_N(u_i))\\
& = & \langle u_i^*[u_i, P_{\Omega_M}\otimes1](\Omega_M\otimes\Omega_N),(\Omega_M\otimes\Omega_N)\rangle\\
& \leq & \|[u_i, P_{\Omega_M}\otimes1](\Omega_M\otimes\Omega_N)\|.
\end{eqnarray*}
As we know that $\|[u_i, P_{\Omega_M} \otimes1](\Omega_M\otimes\Omega_N)\| \to 0$, the proof is complete.
\end{proof}
Before proceeding, we want to remind a few things on basic constructions. Let $B \subset N$ be an inclusion of type ${\rm II_1}$ factors; we denote by $e_B : L^2(N) \to L^2(B)$ the orthogonal projection, and $E_B : N \to B$ the unique trace-preserving conditional expectation from $N$ onto $B$. If we denote by $J_N$ the canonical anti-unitary on $L^2(N)$ associated with the trace $\tau_N$, it is well known that 
\begin{eqnarray*}
\left\langle N, e_B \right\rangle & = & (N \cup e_B)''\\
& = & J_N\{B' \cap B(L^2(N))\}J_N.
\end{eqnarray*}
Moreover, as $e_B$ commutes with $B$ and as for any $x \in N$, $e_Bxe_B = E_B(x)e_B$, $Ne_BN$ turns out to be a $\ast-$algebra; as the central support of $e_B$ is $1$, $Ne_BN$ is weakly dense in $\left\langle N, e_B \right\rangle$. Then we can define a semifinite canonical trace $\Phi$ on $\left\langle N, e_B \right\rangle$ in the following way:
\begin{equation*}
\forall x, y \in N  \;  \Phi(xe_By) = \tau_N(xy).
\end{equation*}
According to \cite{OP}, we remind that if a type ${\rm II_1}$ factor $M$ is decomposed as $M = M_1 \otimes M_2$ for some type ${\rm II_1}$ factors $M_1$, $M_2$ and $t > 0$ then, modulo unitary conjugacy, there exists a unique decomposition $M = M_1^t \otimes M_2^{1/t}$, such that $p_1M_1p_1 \vee p_2M_2p_2$ and $q_1M_1^tq_1 \vee q_2M_2^{1/t}q_2$ are unitary conjugate in $M$ for any projections $p_i \in \Proj(M_i)$, $i = 1, 2$ and $q_1 \in \Proj(M_1^t)$, $q_2 \in \Proj(M_2^{1/t})$, with $\tau(p_1)/\tau(q_1) = \tau(q_2)/\tau(p_2) = t$.

We are now able to state the result by Popa of conjugation of subfactors in type ${\rm II_1}$ factors mentioned in the introduction (see \cite{OP} for other results of the same kind). Our complete classification is based on the following theorem. We wish to gratefully thank Sorin Popa for allowing us  to present it here. We should mention that Sorin Popa has recently given a proof of his result in \cite{popaduff}. For the sake of completeness, we shall give a proof of this theorem.
\begin{theo} \label{popa}
Let $R_1$ and $R_2$ be two copies of the hyperfinite type ${\rm II_1}$ factor; let $N_1$ and $N_2$ be two full factors of type ${\rm II_1}$. Let us assume that the factor $N$ has both of the following decompositions: $N = R_1 \otimes N_1 = R_2 \otimes N_2$. Then there exist $t > 0$ and a unitary $u \in N$ such that $R_2 = u^*(R_1)^{1/t}u$ and $N_2 = u^*(N_1)^tu$. 
\end{theo}
\begin{proof}
First of all, we are looking at the inclusion $N_1 \subset R_1 \otimes N_1 = N$; we know now, according to Lemma $\ref{inclusion}$, that for any central sequence of unitaries $(u_i)$ in $N$, 
\begin{equation*}
\sup_{x \in N_1, \: \|x\| \leq 1}\|(x\otimes1)u_i - u_i(x\otimes1)\|_2 \to 0.
\end{equation*}
If we write the hyperfinite type ${\rm II_1}$ factor $R_2 = \bigotimes_{n \geq 1}(\Mat_2(\C), \tau)$, and if we denote by $R_2^{(n_0)}$ the product $\bigotimes_{n \geq n_0}(\Mat_2(\C), \tau)$, we get immediatly that for any $n_0 \geq 1$, $R_2 = \Mat_{2^{n_0}}(\C) \otimes R_2^{(n_0)}$. Therefore, as $N = R_2 \otimes N_2$, we obtain that
\begin{equation*}
\forall \varepsilon > 0 \; \exists n_0 \in \N \; \forall u \in U(N_1) \; \forall v \in U(R_2^{(n_0)}) \; \|uv - vu\|_2 < \varepsilon.
\end{equation*}
Let $\varepsilon = \frac{1}{2}$. So, we know there exists $n_0 \in \N$ such that for any $u \in U(N_1)$ and any $v \in U(R_2^{(n_0)})$, $\|uv - vu\|_2 < \frac{1}{2}$. For an inclusion $B \subset N$ of type ${\rm II_1}$, we know according to Proposition $1.3.2$ in \cite{popa2003}, that for any $x \in N$
\begin{equation*}
E_{B' \cap N}(x) \in \overline{co}^w\{u^*xu, u \in U(B)\}.
\end{equation*}
Here, let $B = N_1$, then $B' \cap N = R_1$. Noticing that  for any $u \in U(N_1)$ and any $v \in U(R_2^{(n_0)})$, $\|uv - vu\|_2 = \|v - u^*vu\|_2$, we get that for any $v \in U(R_2^{(n_0)})$, $\|v - E_{R_1}(v)\|_2 \leq \frac{1}{2}$ (the same technique is used in \cite{popa2001}). Let now $a \in \overline{co}^w\{v^*e_{R_1}v, v \in U(R_2^{(n_0)})\}$ of minimal $\|\cdot\|_{2, \Phi}-$norm. We are going to prove that $a$ is not $0$. From the one hand, for any $u \in N$, the basic construction for $R_1 \subset N$ gives us $e_{R_1}u^*e_{R_1}ue_{R_1} = E_{R_1}(u^*)E_{R_1}(u)e_{R_1}$. From the other hand, as for any $v \in U(R_2^{(n_0)})$, $\|v - E_{R_1}(v)\|_2 \leq \frac{1}{2}$, we get for any $v \in U(R_2^{(n_0)})$, $\|E_{R_1}(v)\|_2 \geq \frac{1}{2}$. Thus, for any $v \in U(R_2^{(n_0)})$,
\begin{eqnarray*}
\Phi(e_{R_1}v^*e_{R_1}ve_{R_1}) & = & \Phi(E_{R_1}(v^*)E_{R_1}(v)e_{R_1})\\
& = & \tau_N(E_{R_1}(v^*)E_{R_1}(v))\\
& = & \|E_{R_1}(v)\|_2^2\\
& \geq & \frac{1}{4},
\end{eqnarray*}
therefore $\Phi(e_{R_1}ae_{R_1}) \geq \frac{1}{4}$, and $a \neq 0$. Furthermore, as $a$ is of minimal $\|\cdot\|_{2, \Phi}-$norm, necessarily we have for any $v \in U(R_2^{(n_0)})$, $a = v^*av$. Finally, we have found $a \in \langle N, e_{R_1}\rangle^+\cap(R_2^{(n_0)})'$ with $a \neq 0$ and $\Phi(a) < \infty$. We can notice that $\langle N, e_{R_1}\rangle = R_1 \otimes B(L^2(N_1))$; so, if we denote by $e$ the spectral projection of $a$ corresponding to the interval $[\|a\|/2, \|a\|]$, then $e$ is a finite non-zero projection in $R_1 \otimes B(L^2(N_1))$, which commutes with $R_2^{(n_0)}$. We can proceed by using the proof of Proposition $12$ of \cite{OP}. The Hilbert space $\h = eL^2(N)$ is an $R_2^{(n_0)} - R_1$ Hilbert bimodule and $\dim \h_{R_1} < \infty$. As $R_1' \cap N = N_1$ is a factor, Proposition $12$ of \cite{OP} claims that there exists $s > 0$ and $u \in U(N)$ such that $uR_2^{(n_0)}u^* \subset R_1^s$. Thus, up to a conjugation by a unitary in $N$, we get $R_2^{(n_0)} \subset R_1^s$. But by construction, we know that $R_2 = \Mat_{2^{n_0}}(\C) \otimes R_2^{(n_0)}$, so up to a conjugation by a unitary in $N$, we get $R_2^{1/2^{n_0}} \subset R_1^s$, i.e. up to a conjugation by a unitary in $N$, we have $R_2 \subset R_1^{s'}$ with $s' = 2^{n_0}s$. So, we get that $N_2 = R_0 \vee N_1^{1/s'}$ with $R_0 = R_2' \cap R_1^{s'}$. Since $N_2 = R_0 \otimes N_1^{1/s'}$, $R_0$ is a subfactor of $R_1^{s'}$. As $N_2$ is a full factor, $R_0$ has to be finite-dimensional. If $d^2 = \dim R_0$, $t = d/s'$, and if we decompose $N$ as $R_1^{1/t} \otimes N_1^t$, we obtain up to a conjugation by a unitary in $N$ that $N_2 = N_1^{t}$ and $R_2 = R_1^{1/t}$.
\end{proof}
The following lemma is an easy consequence of the previous theorem. It is this result we will use in the proof of Theorem $\ref{bijection}$.
\begin{lem}\label{iso}
Let $R_1^{\infty}$ and $R_2^{\infty}$ be two copies of the hyperfinite type ${\rm II_{\infty}}$ factor; let $N_1^{\infty}$ and $N_2^{\infty}$ be two full factors of type ${\rm II_{\infty}}$. Let $\a : R_1^{\infty} \otimes N_1^{\infty} \to R_2^{\infty} \otimes N_2^{\infty}$ be an isomorphism. Then, there exist a unitary $u \in R_2^{\infty} \otimes N_2^{\infty}$ and two isomorphisms $\b : R_1^{\infty} \to R_2^{\infty}$, $\g : N_1^{\infty} \to N_2^{\infty}$  such that for any $x \in R_1^{\infty} \otimes N_1^{\infty}$, $u^*\a(x)u = (\b \otimes \g)(x)$.
\end{lem}
\begin{proof}
Let $\a$ be an isomorphism from $R_1^{\infty} \otimes N_1^{\infty}$ onto $R_2^{\infty} \otimes N_2^{\infty}$. Let $e_1$ and $f_1$ be two finite projections in $R_1^{\infty}$ et $N_1^{\infty}$. It is clear that $\a(e_1 \otimes f_1)$ is a finite projection in $R_2^{\infty} \otimes N_2^{\infty}$; thus there exist two finite projections $e_2$ and $f_2$ in $R_2^{\infty}$ and $N_2^{\infty}$ such that $\a(e_1 \otimes f_1) \sim e_2 \otimes f_2$. If we denote by $p_1$ the projection $\a(e_1 \otimes f_1)$, as $R_2^{\infty} \otimes N_2^{\infty}$ is a properly infinite factor, there exists a family $(p_n)$ of pairwise orthogonal and mutually equivalent projections such that $\sum p_n = 1$. If we do again the same thing with $q_1 = e_2 \otimes f_2$, we obtain another family $(q_n)$ of pairwise orthogonal and mutually equivalent projections such that $\sum q_n = 1$. But, as $p_1 \sim q_1$, we obtain for any $n \in \N$, $p_n \sim q_n$. Let $u_n \in R_2^{\infty} \otimes N_2^{\infty}$ be the partial isometry such that $u_n^*u_n = p_n$ and $u_nu_n^* = q_n$; then $u = \sum u_n$ is a unitary in $R_2^{\infty} \otimes N_2^{\infty}$ and $up_1u^* = q_1$. Consequently, up to a conjugation by a unitary, we can assume that $\a(e_1 \otimes f_1) = e_2 \otimes f_2$.

From now on, we denote by $R_1$ the factor $(R_1^{\infty})_{e_1}$ and $N_1$ the factor $(N_1^{\infty})_{f_1}$. We obtain immediately that $\a(R_1 \otimes N_1) = (R_2^{\infty})_{e_2} \otimes (N_2^{\infty})_{f_2}$, because $\a(e_1 \otimes f_1) = e_2 \otimes f_2$. According to Theorem $\ref{popa}$, as $\a(R_1 \otimes f_1)$, $(R_2^{\infty})_{e_2}$ are two copies of the hyperfinite factor of type ${\rm II_1}$, and $a(e_1 \otimes N_1)$, $(N_2^{\infty})_{f_2}$ are two full factors of type ${\rm II_1}$, we know there exists a unitary $u \in \a(R_1 \otimes N_1)$ and $t > 0$ such that:
\begin{eqnarray*}
u\a(R_1 \otimes f_1)u^* & = & \left\{(R_2^{\infty})_{e_2}\right\}^{1/t}\\
\and u\a(e_1 \otimes N_1)u^* &  = & \left\{(N_2^{\infty})_{f_2}\right\}^{t}.
\end{eqnarray*}
Moreover, there exist a projection $e'_1$ in $R_2^{\infty}$ and a projection $f'_1$ in $N_2^{\infty}$ such that:
$$\left\{(R_2^{\infty})_{e_2}\right\}^{1/t} = (R_2^{\infty})_{e'_1} \and \left\{(N_2^{\infty})_{f_2}\right\}^{t} = (N_2^{\infty})_{f'_1}.$$
Consequently, we have proved:
\begin{eqnarray*}
u\a(R_1 \otimes f_1)u^* & = & (R_2^{\infty})_{e'_1} \otimes f'_1\\
\and u\a(e_1 \otimes N_1)u^* &  = & e'_1 \otimes (N_2^{\infty})_{f'_1}.
\end{eqnarray*}
Thus, up to a conjugation by a unitary, we can assume that:
\begin{eqnarray*}
\a(R_1 \otimes f_1) & = & (R_2^{\infty})_{e'_1} \otimes f'_1\\
\and \a(e_1 \otimes N_1) &  = & e'_1 \otimes (N_2^{\infty})_{f'_1}.
\end{eqnarray*}
Let $(e_{ij})_{i, j \geq 1}$ be a system of matrix unit in $R_1^{\infty}$ such that $e_{11} = e_1$, $R_1^{\infty} = V^*\left(R_1 \otimes B(l^2)\right)V$ with the unitary $V$ defined in the following way (we shall assume that $R_1^{\infty}$ acts on the Hilbert space $\h$):
$$V \: : \: \fonction{\h}{e_1\h \otimes l^2}{h}{(e_{1n}h)}.$$
In the same way, we denote by $(f_{ij})_{i, j \geq 1}$, $(e'_{ij})_{i, j \geq 1}$, $(f'_{ij})_{i, j \geq 1}$ the systems of matrix unit associated respectively with $N_1^{\infty}$, $R_2^{\infty}$, $N_2^{\infty}$. Let us denote by $U$ the unitary $\sum_{i, j}(e'_{i1}\otimes f'_{j1})\a(e_{1i} \otimes f_{1j})$. The computation of $U\a(e_{ik} \otimes f_{jl})U^*$ gives us:
\begin{eqnarray*}
U\a(e_{ik} \otimes f_{jl})U^* & = & (e'_{i1} \otimes f'_{j1})\a(e_{1i} \otimes f_{1j})\a(e_{ik} \otimes f_{jl})\a(e_{k1} \otimes f_{l1})(e'_{1k} \otimes f'_{1l})\\
& = & (e'_{i1} \otimes f'_{j1})\a(e_{11} \otimes f_{11})(e'_{1k} \otimes f'_{1l})\\
& = & (e'_{i1} \otimes f'_{j1})(e'_{11} \otimes f'_{11})(e'_{1k} \otimes f'_{1l})\\
& = & (e'_{ik} \otimes f'_{jl}).
\end{eqnarray*}
Let $\b : R_1^{\infty} \to R_2^{\infty}$ defined in the following way: $\b(x) = \a(x)$ for all $x \in R_1$ and $\b(e_{ij}) = e'_{ij}$ for all $i, j \geq 1$. It is clear that $\b$ is an isomorphism from $R_1^{\infty}$ onto $R_2^{\infty}$. In the same way, we define $\g : N_1^{\infty} \to N_2^{\infty}$ by: $\g(x) = \a(x)$ for all $x \in N_1$ and $\g(f_{ij}) = f'_{ij}$ for all $i, j \geq 1$. Once again, $\g$ is an isomorphism from $N_1^{\infty}$ onto $N_2^{\infty}$. Finally, with these notations, for any $x \in R_1^{\infty} \otimes N_1^{\infty}$, we have:
$$U\a(x)U^* = (\b \otimes \g)(x).$$
\end{proof}
\begin{proof}[Proof of Theorem $\ref{bijection}$]
We keep the notations introduced at the beginning of this section, $C = N \indice{\Join}_{\g} L^{\infty}(X)$ and $C' = N' \indicebis{\Join}_{\g'} L^{\infty}(X')$. 
We have seen, according to Theorem $\ref{construction}$, that the core of $C$ is isomorphic to $\left(L^{\infty}(X) \rtimes_{\g} \R\right) \otimes N$ and the dual action $\theta_s$ is given by $\theta_s = \widehat{\g}_{-s} \otimes \a_s$. The same thing is true for $C'$. As $C \simeq C'$, $(\Core(C), \theta)$ and $(\Core(C'), \theta')$ are weakly equivalent, according to Corollary $8.4$ of \cite{takesaki73}. That means that there exist an isomorphism $\pi : \Core(C) \to \Core(C')$ and a strongly continuous family of unitaries $(u_s)$ of $\Core(C)$ such that $(u_s)$ is a one-cocycle for $\t_s$ and for any $x \in \Core(C)$ and any $s \in \R$,
$$\pi^{-1}\theta'_s\pi(x) = u_s\theta_s(x)u_s^*.$$
But, according to Lemma \ref{iso}, we know that there exist a unitary $u \in \Core(C')$ and two isomorphisms $\pi_1 : L^{\infty}(X) \rtimes_{\g} \R \to L^{\infty}(X') \rtimes_{\g'} \R$, $\pi_2 : N \to N'$ such that for any $x \in \Core(C)$,
$$u\pi(x)u^* = (\pi_1 \otimes \pi_2)(x).$$
Thus we obtain, for any $x \in \Core(C)$ any $s \in \R$:
\begin{eqnarray*}
\theta'_s(\pi(x)) & = & \pi(u_s)\pi(\theta_s(x))\pi(u_s)^*\\
\theta'_s(u)^*\theta'_s((\pi_1 \otimes \pi_2)(x))\theta'_s(u) & = & \pi(u_s)u^*(\pi_1 \otimes \pi_2)(\theta_s(x))u\pi(u_s)^*.
\end{eqnarray*}
If $v_s = (\pi_1 \otimes \pi_2)^{-1}(\theta'_s(u)\pi(u_s)u^*)$, for any $x \in \Core(C)$ and any $s \in \R$, we have
$$(\pi_1 \otimes \pi_2)^{-1}\theta'_s((\pi_1 \otimes \pi_2)(x)) = v_s\theta_s(x)v_s^*.$$
If we take now $x = z \otimes 1$, for $z \in L^{\infty}(X) \rtimes_{\g} \R$, we have
$$(\pi_1)^{-1}\left(\widehat{\g'}_{-s}(\pi_1(z))\right) \otimes 1 = v_s(\widehat{\g}_{-s}(z) \otimes 1)v_s^*.$$
Let us denote $\b_s = (\pi_1)^{-1}\widehat{\g'}_{-s}\pi_1\widehat{\g}_{s}$; it is an automorphism of $L^{\infty}(X) \rtimes_{\g} \R$. Moreover, for any $y \in L^{\infty}(X) \rtimes_{\g} \R$, we have
$$(\b_s(y) \otimes 1)v_s = v_s(y \otimes 1).$$
Now, we are going to use a classical technique. Let $\phi \in N_*$; we have, for any $y \in L^{\infty}(X) \rtimes_{\g} \R$, the equality: 
$$\b_s(y)(\id \otimes \phi)(v_s) = (\id \otimes \phi)(v_s)y.$$
As $v_s \neq 0$, there exists necessarily $\phi \in N_*$ such that $(\id \otimes \phi)(v_s) \neq 0$. A classical lemma allows us to claim that $\b_s$ is an inner automorphism of $L^{\infty}(X) \rtimes_{\g} \R$. Therefore, there exists a family of unitaries $(w_s)$ in $L^{\infty}(X) \rtimes_{\g} \R$ (we can choose it strongly continuous because $s \mapsto \b_s$ is a continuous map), such that for any $z \in L^{\infty}(X) \rtimes_{\g} \R$,
\begin{equation}\label{cocycle}
(\pi_1)^{-1}\widehat{\g'}_s\pi_1(z) = w_s\widehat{\g}_s(z)w_s^*.
\end{equation}
But, we have to be careful here because $(w_s)$ is not in general a one-cocycle for $\widehat{\g}$; for $s$ and $t \in \R$, 
\begin{eqnarray*}
(\pi_1)^{-1}\widehat{\g'}_s\widehat{\g'}_t\pi_1 & = & \Ad(w_s\widehat{\g}_s(w_t))\widehat{\g}_{s + t}\\
\and (\pi_1)^{-1}\widehat{\g'}_{s + t}\pi_1 & = & \Ad(w_{s + t})\widehat{\g}_{s + t}.
\end{eqnarray*}
Consequently, $\Ad(w_s\widehat{\g}_s(w_t))\widehat{\g}_{s + t} =  \Ad(w_{s + t})\widehat{\g}_{s + t}$, and we obtain that $w_{s + t}^*w_s\widehat{\g_s}(w_t)$ is in the center of $L^{\infty}(X) \rtimes_{\g} \R$; As $L^{\infty}(X) \rtimes_{\g} \R$ is a factor, we have finally proved that $w_{s + t}^*w_s\widehat{\g_s}(w_t) \in \mathbf{T}$. Hence, there exists $\omega(s, t) \in \mathbf{T}$, such that $w_{s + t} = \omega(s, t)w_s\widehat{\g_s}(w_t)$. Moreover, we can easily see that $\omega(s, t)$ satisfies a 2-cocycle relation. But, it is well-known that $Z^2(\R, \mathbf{T}) = B^2(\R, \mathbf{T})$, therefore $\omega(s, t)$ is a coboundary and there exists a map $\la : \R \to \mathbf{T}$, such that for any $s$, $t \in \R$, $\omega(s, t) = \overline{\la(s + t)}\la(s)\la(t)$. We can notice that if we multiply $w_s$ by $\la(s)$, $\la(s)w_s$ becomes a one-cocycle for $\widehat{\g}$ and we do not change the Equation (\ref{cocycle}). Finally, up to a multiplication by a scalar $\la(s) \in \mathbf{T}$, we can assume that $w_s$ is a one-cocycle for $\widehat{\g}$.

We can apply now Proposition $4.2$ of \cite{takesaki73}. We obtain the existence of an isomorphism $\tilde{\pi} : L^{\infty}(X) \otimes B(L^2(\R)) \to L^{\infty}(X') \otimes B(L^2(\R))$ which intertwines both of the bidual actions $\widehat{\widehat{\g}}$ and $\widehat{\widehat{\g'}}$, i.e. for any $x \in L^{\infty}(X) \otimes B(L^2(\R))$,
$$\tilde{\pi}\widehat{\widehat{\g}}_s(x) = \widehat{\widehat{\g'}}_s\tilde{\pi}(x).$$
Moreover we have both of the equalities $\widehat{\widehat{\g}}_s = \g_s \otimes \Ad(\la_s^*)$ and $\widehat{\widehat{\g'}}_s = \g'_s \otimes \Ad(\la_s^*)$. We can notice that $L^{\infty}(X) = \centre(L^{\infty}(X) \otimes B(L^2(\R)))$ and $L^{\infty}(X') = \centre(L^{\infty}(X') \otimes B(L^2(\R)))$. Let us denote by $\rho : L^{\infty}(X) \to L^{\infty}(X')$ the isomorphism obtained from $\tilde{\pi}$ by taking restrictions to the centers. Finally, we obtain for any $s \in \R$ and $x \in L^{\infty}(X)$,
$$\rho\g_s(x) = \g'_s\rho(x).$$
Therefore, we have proved that the flows $(X, \g)$ and $(X',\g'_t)$ are conjugate; we can prove exactly in the same way that $(N, \a)$ and $(N', \a')$ are weakly equivalent.
\end{proof}

\section{Computation of Connes' $\tau$ Invariant}

For this section, we are going to keep the same notations as before. We shall denote by $C$ the factor of type ${\rm III_1}$, $N \indice{\Join}_{\g} L^{\infty}(X)$; we remind that $C \otimes B(L^2(\R))$ is canonically isomorphic to $(L^{\infty}(X) \otimes P) \rtimes_{\g \otimes \s^{\omega}} \R$ with $P = N \rtimes_{\a} \R$. We assume that the flow $\g$ is free and ergodic and finite or infinite measure-preserving; the action $\g \otimes \s^{\omega}$ will be still denoted by $\b$ and $\la_s$ are the unitaries which implement the action of $\R$ on $L^{\infty}(X) \otimes P$. Furthermore, we shall assume that $P$ is a free Araki-Woods factor \cite{shlya97}; but $N$ is not necessarily full. At last, we shall denote by $\s$ the canonical embedding of $(L^{\infty}(X) \otimes P) \rtimes_{\b} \R$ into $M = P \otimes L^{\infty}(X) \otimes B(L^2(\R))$. First of all, let us remind the definition of Connes' $\tau$ invariant.
\begin{df}
Let $C$ be a factor of type ${\rm III_1}$. Let $\varepsilon : \Aut(C) \to \Out(C) = \Aut(C)/\Inn(C)$ be the canonical projection. Let $\psi$ be any weight on $C$. The mapping $\delta : \R \to \Out(C)$, with $\delta(t) = \varepsilon(\s_t^{\psi})$ is independent of the choice of the weight $\psi$ thanks to Theorem $1.2.1$ of \cite{connes73}. The $\tau$ invariant of $C$, denoted by $\tau(C)$ is the \emph{weakest topology} on $\R$ that makes the map $\delta$ continuous.
 \end{df}
Although $C$ is not a full factor, we can nevertheless give an explicit computation of Connes' $\tau$ invariant in many cases. We want to remind a few things about free Araki-Woods factors \cite{shlya97}. Let $\h_{\R}$ be a real Hilbert space and $(U_t)$ be an orthogonal representation of $\R$ on $\h_{\R}$. Let $\h = \h_{\R} \otimes_{\R} \C$ be the complexified Hilbert space. If $A$ is the infinitesimal generator of $(U_t)$ on $\h$, we remind that $j : \h_{\R} \to \h$ defined by $j(\zeta) = (\frac{2}{A^{-1} + 1})^{1/2}\zeta$ is an isometric embedding of $\h_{\R}$ into $\h$. Let $K_{\R} = j(\h_{\R})$. For $\xi \in \h$, we denote by $l(\xi)$ the creation operators on the Fock space $\mathcal{F}(\h)$ and $s(\xi)$ their real part. By definition, $\Gamma(\h_{\R}, U_t)'' = \{s(\xi), \xi \in K_{\R}\}''$. 
\begin{df}
Let $\h_{\R}$ be a real Hilbert space and $(U_t)$ be an orthogonal representation of $\R$ on $\h_{\R}$. Let $P = \Gamma(\h_{\R}, U_t)''$ be the free Araki-Woods factor associated with $(U_t)$ and $\h_{\R}$  \cite{shlya97}. We shall say that $P$ satisfies the \emph{condition (M)} of mixing if there exist $\xi$, $\eta \in \h_{\R}$, such that the continuous function $f$ defined by $f(t) = \langle U_tj\xi, j\eta\rangle$ vanishes at infinity and $f \neq 0$.
\end{df}
For example, if $\h_{\R} = L^2(\R, \R)$ and $U_t = \la_t$ for any $t \in \R$, the factor $\Gamma(L^2(\R, \R), \la_t)''$ satisfies the condition (M). We remind that $\Gamma(L^2(\R, \R), \la_t)''$ is nothing but $N \rtimes_{\a} \R$, with $N = L(\mathbf{F}_{\infty}) \otimes \bh$ and $\a$ the trace-scaling automorphism group of R\u{a}dulescu \cite{{radulescu1991}, {shlya98}}. The aim of this section is to prove the following theorem: 
\begin{theo}\label{computation}
Let $P = \Gamma(\h_{\R}, U_t)''$ be a free Araki-Woods factor satisfying the condition $(M)$. For $C = N\indice{\Join}_{\g}L^{\infty}(X)$,
$\tau(C)$ = the weakest topology on $\R$ that makes the map $t \mapsto \g_t$ continuous for the $u$-topology on $\Aut(L^{\infty}(X))$.
\end{theo} 
The proof of this result is based on a kind of \textquotedblleft $14\varepsilon$ lemma\textquotedblright, it will be stated in the Appendix (Lemma $\ref{14epsilon}$).

From now on, we are going to use notations and results of \cite{connes74}. For $C$ our factor of type ${\rm III_1}$, for any normal faithful state $\varphi$ on $C$ and for any free ultrafilter $\omega$ on $\N$, we shall denote by $A_{\varphi, \omega}$ the norm closed $\ast-$subalgebra of $l^{\infty}(\N, C)$ of all sequences $(x_n)_{n \in \N}$ such that $\|[x_n, \varphi]\| \to 0$ when $n \to \omega$. Let $I_{\omega} =  \{(x_n)_{n \in \N}, x_n \to 0 \:*-\mbox{strongly when } n \to \omega\}$. It is a two-sided ideal of $l^{\infty}(\N, C)$. Moreover, $A_{\varphi, \omega} \cap I_{\omega}$ is a two-sided ideal of the $C^*-$algebra $A_{\varphi, \omega}$, and the canonical quotient map will be denoted by $\rho_{\varphi, \omega}$. We shall denote by $C^{\omega}$ the \emph{ultraproduct} of $C$ along $\omega$, i.e. the quotient of $l^{\infty}(\N, C)$ by the two-sided ideal $I_{\omega}$; we shall denote by $C_{\varphi, \omega}$ the quotient of the $C^*-$algebra $A_{\varphi, \omega}$ by the two-sided ideal $A_{\varphi, \omega} \cap I_{\omega}$. We know, according to Proposition $2.2$ of \cite{connes74}, that $C_{\varphi, \omega}$ is a finite von Neumann algebra. Let $\varphi$ be a given faithful normal state on $C$ and $\mathcal{D}$ the set of faithful normal states on $C$ with $\a\varphi \leq \psi \leq \a^{-1}\varphi$ for some $\a > 0$. Let
\begin{equation*}
C_{\omega} = \bigcap_{\psi \in \mathcal{D}}\rho_{\varphi, \omega} (A_{\varphi, \omega} \cap A_{\psi, \omega}).
\end{equation*}
According to Theorem $2.9$ of \cite{connes74}, $C_{\omega}$ is a finite von Neumann called the \emph{asymptotic centralizer} of $C$ at $\omega$. Before proceeding, we want to remind a classical result on compact operators on the Hilbert space $L^2(\R)$.
\begin{prop}\label{compact}
Let $A \subset B(L^2(\R))$ be a unital $C$*-algebra. Assume that for any $t \in \R$, $\la_t \in A$ and that there exists $f \in C_o(\R)$, $f \neq 0$, such that $M_f \in A$ ($M_f$ is the multiplication operator). Then $K(L^2(\R)) \subset A$.
\end{prop}
\begin{proof}
We identify $g$ with $M_g$ for any $g \in C_o(\R)$. Since $f \neq 0$ and all its translations belong to $A \cap C_o(\R)$, $A \cap C_o(\R)$ is a sub$-C^*-$algebra of $C_o(\R)$ which separates the points. Consequently, $A \cap C_o(\R) = C_o(\R)$ and thus $C_o(\R)\subset A$. Moreover, $C^*_r(\R) \subset A$ because for any $t \in \R$, $\la_t \in A$. Therefore, $K(L^2(\R)) = [C^*_r(\R)C_o(\R)] \subset A$.
\end{proof}
For any subset $E \subset C$ and a sequence $(x_k)$ of elements in $C$, we shall say that $(x_k)$ \emph{almost commutes} with $E$ if for any $a \in E$, $[x_k, a] \to 0$ $\ast-$strongly. Let $C = N \indice{\Join}_{\g} L^{\infty}(X)$; we know that $C \otimes B(L^2(\R))$ is canonically isomorphic to $(L^{\infty}(X) \otimes P) \rtimes_{\b} \R$. We remind that we have denoted by  $\tilde{\phi}$ the dual weight on $(L^{\infty}(X) \otimes P) \rtimes_{\b} \R$ of the weight $\phi = \tau \otimes \omega$ on $L^{\infty}(X) \otimes P$. Let $\la_t$ be the unitaries in $(L^{\infty}(X) \otimes P) \rtimes_{\b} \R$ which implement the action $\b = \g \otimes \s^{\omega}$ of $\R$ on $L^{\infty}(X) \otimes P$. As the family $(\la_t)$ is in the centralizer of the weight $\tilde{\phi}$, it is a one-cocycle for the action $\s^{\tilde{\phi}}$. We know according to Theorem $1.2.4$ of \cite{connes73} that there exists a faithful, normal, semifinite weight $\psi$ on $(L^{\infty}(X) \otimes P) \rtimes_{\b} \R$ such that:
$$\forall t \in \R, \: \s_t^{\psi} = \la_t^*\s_t^{\tilde{\phi}}\la_t.$$
We remind that we have shown in the first section that for any $t \in \R$ and any $x \in L^{\infty}(X) \otimes P$,
\begin{equation}\label{poids}
\s_t^{\psi}(\pi_{\b}(x))  =  \pi_{\b}((\g_{-t} \otimes \id)(x)).
\end{equation}
Thanks to Equation $(\ref{poids})$, we know that $\Cent(\psi) = (\pi_{\s^{\omega}}(P)\cup\la(\R))'' \simeq P\rtimes_{\s^{\omega}}\R$.
We are proving now, thanks to Proposition $\ref{conclusion}$, the following technical lemma which will turn out to be essential in proof of Theorem~$\ref{computation}$:
\begin{lem}\label{technique}
Let $C$ be as in Theorem $\ref{computation}$, and $\psi$ the weight on $(L^{\infty}(X) \otimes P) \rtimes_{\b} \R$ which satisfies Equation $(\ref{poids})$. Let $\omega$ be any free ultrafilter on $\N$. Any bounded sequence $(u_n)$ of $(L^{\infty}(X) \otimes P) \rtimes_{\b} \R$ which almost commutes with $\Cent(\psi)$ is centralising. In a shorter way, we have the following equality:
$$((L^{\infty}(X) \otimes P) \rtimes_{\b} \R)^{\omega}\cap\Cent(\psi)' = ((L^{\infty}(X) \otimes P) \rtimes_{\b} \R)_{\omega}.$$
\end{lem}
\begin{proof}
The inclusion $\supset$ is trivial. We have just to prove $\subset$. We assume as in Theorem $\ref{computation}$, that $P = \Gamma(\h_{\R}, U_t)''$ is the free Araki-Woods factor associated with the orthogonal representation $(U_t)$. Let $\Omega$ be the vacuum vector and let $\varphi_U = \langle \cdot \Omega, \Omega\rangle$ be the free quasi-free state on $P$, and $(\s_t)$ the modular group associated with $\varphi_U$. We remind that we have for any $t \in \R$, and any $\zeta \in K_{\R}$, $\s_t(s(\zeta)) = s(U_t\zeta)$. We remind that we take $M = P \otimes L^{\infty}(X) \otimes B(L^2(\R))$ and the embedding of $(L^{\infty}(X) \otimes P) \rtimes_{\b} \R$ into $M$ is denoted by $\s$. Let $(u_n)$ be a sequence of elements in $(L^{\infty}(X) \otimes P) \rtimes_{\b} \R$ which almost commutes with $\Cent(\psi)$ and such that $\|u_n\| \leq 1$ for any $n \in \N$. Thus $(\s(u_n))$ almost commutes with $\s(\Cent(\psi))$ and in particular $(\s(u_n))$ almost commutes with $\s(P)$. Hence, according to Proposition \ref{conclusion}, there exists a bounded sequence $(v_n)$ of elements in $L^{\infty}(X)\rtimes_{\g}\R$ such that $\s(u_n) - 1\otimes v_n \to 0$ $\ast-$strongly in $M = P \otimes L^{\infty}(X) \otimes B(L^2(\R))$ and $\|v_n\| \leq 1$ for any $n \in \N$. Actually, $v_n$ is nothing but $E(u_n)$, for any $n \in \N$, with $E : M \to L^{\infty}(X)\otimes B(L^2(\R))$ the conditional expectation $\varphi_U \otimes \Id$. We have somewhat reduced the difficulty of the problem: instead of looking at the sequence $(u_n)$ in $(L^{\infty}(X) \otimes P) \rtimes_{\b} \R$, we are looking at the sequence $(v_n)$ in $L^{\infty}(X)\rtimes_{\g}\R$. But we want to go further; we want to prove that there exists a sequence $(w_n)$ in $L^{\infty}(X)$ such that $\s(u_n) - 1\otimes w_n\otimes1 \to 0$ $\ast-$strongly in $M$. If we are able to prove such a result, we are done. Indeed, as $1\otimes w_n\otimes1 \in \centre(M)$ for any $n \in \N$, the sequence $(1\otimes w_n\otimes1)$ turns out to be centralising in $M$; but as $\s(u_n) - 1\otimes w_n\otimes1 \to 0$ $\ast-$strongly in $M$, $(\s(u_n))$ is also centralising in $M$. Therefore, $(u_n)$ is centralising in $(L^{\infty}(X) \otimes P) \rtimes_{\b} \R$.

As $\s(u_n) - 1\otimes v_n \to 0$ $\ast-$strongly in $M$ and $(u_n)$ almost commutes with $\Cent(\psi)$, for any $y \in \Cent(\psi)$, $[1\otimes v_n, y_{13}] \to 0$ $\ast-$strongly in $M$ (we are regarding $y$ as a element in $M = P \otimes L^{\infty}(X) \otimes B(L^2(\R))$, that is why it is denoted by $y_{13}$). Thus for any normal state $\phi$ on $P$, $[v_n, 1\otimes (\phi \otimes \id)(y)] \to 0$ $\ast-$ strongly in $L^{\infty}(X) \otimes B(L^2(\R))$. Let $B = C^*((\phi\otimes\Id)(y), y \in P\rtimes_{\s^{\omega}}\R, \phi \in P_*)$. We see that $B \subset B(L^2(\R))$ is a unital $C^*-$algebra, and for any $t \in \R$, $\la_t \in B$. 
Now, since $P$ satisfies condition (M), we know that there exist $\xi$, $\eta \in \h_{\R}$ such that the continuous function $f$ defined by $f(t) = \langle U_t\xi, \eta\rangle$ vanishes at infinity and $f \neq 0$. Let $\xi' = j(\xi)$ and $\eta' = j(\eta)$. Let $\phi = \varphi_U(s(\eta')\cdot) \in P_*$ and $x = (\s_t(s(\xi')))_{t \in \R} \in \Core(P)$, with $\Core(P) = P\rtimes_{\s^{\omega}}\R$. We get immediatly that for any $t \in \R$, 
\begin{eqnarray*}
(\phi\otimes\Id)(x)(t) & = & \langle s(\eta')\s_t(s(\xi'))\Omega, \Omega\rangle \\
            & = &   \langle \s_t(s(\xi'))\Omega, s(\eta')\Omega\rangle \\
            & =  &  \langle s(U_t\xi')\Omega, s(\eta')\Omega\rangle \\
            & = & \frac{1}{4} \langle U_t\xi', \eta'\rangle.
\end{eqnarray*}
Consequently, $(\phi \otimes \Id)(x) \in C_o(\R)$  and $(\phi \otimes \Id)(x) \neq 0$. Then, according to Proposition \ref{compact}, $K(L^2(\R)) \subset B$. Hence we have proved that for any $y \in K(L^2(\R))$, $[v_n, 1\otimes y] \to 0$ $\ast-$strongly. Let $\xi \in L^2(\R)$ such that $\|\xi\| = 1$; for any $\mu \in L^2(\R)$ we shall denote again by $\mu$ the map from $\C$ into $L^2(\R)$ such that for any $\la \in \C$, $\mu(\la) = \la\mu$. Thus for $\mu$, $\mu' \in L^2(\R)$, $\t_{\mu', \mu} = \mu' \circ \mu^*$ is a rank-one operator and it is exactly $\langle \cdot, \mu\rangle\mu'$. In particular, these operators are compact. For any $n \in \N$, let $w_n$ be the element of $L^{\infty}(X)$ defined by $(1\otimes\xi^*)v_n(1\otimes\xi)$. For any $\eta \in L^2(X)$, we have 
\begin{eqnarray*}
\|w_n\eta\|^2 & = & \langle w_n^*w_n\eta, \eta\rangle \\
                       & = &  \langle (1\otimes\xi^*)v_n^*(1\otimes\xi)(1\otimes\xi^*)v_n(1\otimes\xi)\eta, \eta\rangle \\
                       & = & \langle v_n^*(1\otimes \t_{\xi, \xi})v_n(\eta \otimes \xi), \eta \otimes \xi\rangle \\
                       & = & \langle v_n^*((1\otimes \t_{\xi, \xi})v_n - v_n(1\otimes \t_{\xi, \xi}))(\eta \otimes \xi), \eta \otimes \xi\rangle + \|v_n(\eta \otimes \xi)\|^2. 
\end{eqnarray*}
Therefore, we get
\begin{equation*}
\left|\|w_n\eta\|^2 - \|v_n(\eta \otimes \xi)\|^2\right| \leq \|(v_n(1\otimes \t_{\xi, \xi}) - (1\otimes \t_{\xi, \xi})v_n)(\eta \otimes \xi)\| \: \|\eta\|.
\end{equation*}
For any $\eta \in L^2(X)$ and $\mu \in L^2(\R)$, we have 
\begin{eqnarray*}
\|(v_n - w_n\otimes1)(\eta \otimes \mu)\| & = & \|(v_n - w_n\otimes1)(1\otimes \t_{\mu, \xi})(\eta \otimes \xi)\| \\
                                                                        & \leq & \|(v_n(1\otimes \t_{\mu, \xi}) - (1\otimes \t_{\mu, \xi})v_n)(\eta \otimes \xi)\|\\
                                                                        &  &  + \|(1\otimes \t_{\mu, \xi})(v_n - w_n\otimes1)(\eta \otimes \xi)\| \\
& \leq & \|(v_n(1\otimes \t_{\mu, \xi}) - (1\otimes \t_{\mu, \xi})v_n)(\eta \otimes \xi)\|\\
&  &  + \|(v_n - w_n\otimes1)(\eta \otimes \xi)\| \: \|\mu\|.
\end{eqnarray*}
But $\|(v_n - w_n\otimes1)(\eta \otimes \xi)\|^2 = \|v_n(\eta \otimes \xi)\|^2 - \|w_n\eta\|^2$. Thus, using the several inequalities we obtained above and the fact that for any $y \in K(L^2(\R))$, $[v_n, 1\otimes y] \to 0$ $\ast-$strongly, for any $\varepsilon > 0$, there exists $n_0 \in \N$ large enough such that for any $n \geq n_0$, 
\begin{equation*}
\|(v_n - w_n\otimes1)(\eta \otimes \mu)\| \leq \varepsilon.
\end{equation*} 
That means exactly that $(v_n - w_n \otimes 1) \to 0$ strongly. We can do the same thing with $v_n^*$ and $w_n^*$ instead of $v_n$ and $w_n$. Hence we have proved that $(v_n - w_n \otimes 1) \to 0$ $\ast-$strongly. The proof is complete.
\end{proof}
\begin{proof}[Proof of Theorem $\ref{computation}$]
Here, as usual, $\Aut(C)$ is endowed with the $u$-topology. Let us denote by $\tau$ the weakest topology on $\R$ that makes the map $t \mapsto \g_t$ continuous for the $u$-topology on $\Aut(L^{\infty}(X))$. According to Equation $(\ref{presqueper})$, as $(\g_{-t})$ and the restriction of $(\s^{\psi}_t)$ on the von Neumann algebra $L^{\infty}(X)$ are exactly the same, it is clear that for any sequence of real numbers $(t_n)$, if $t_n \to 0$ w.r.t. the topology $\tau$ then $t_n \to 0$ then w.r.t. the topology $\tau(C)$. We are going to prove the converse. Let $(t_n)$ be a sequence of real numbers such that $t_n \to 0$ w.r.t. the topology $\tau(C)$. There exists a sequence of unitaries $(u_n)$ in $C$ such that $\Ad(u_n) \circ \s^{\psi}_{t_n} \to \Id$ in $\Aut(C)$. Then the sequence $(u_n)$ almost commutes with $\Cent{\psi}$. So according to Lemma $\ref{technique}$, $(u_n)$ is centralising in $C$. That means exactly that $\Ad(u_n) \to \Id$ in $\Aut(C)$, thus $\s^{\psi}_{t_n} \to \Id$ in $\Aut(C)$. Therefore $t_n \to 0$ w.r.t. the topology $\tau$.
\end{proof}

We wish to end this section by giving a consequence of the previous result. First, we are going to remind a few things about mixing actions on a probability space.
\begin{df}\cite{petersen}
Let $(X, \mu)$ be a probability space, and let $\g$ be a measure-preserving transformation of $(X, \mu)$. Let $\tau$ be the canonical trace on $L^{\infty}(X)$ given by the probability measure $\mu$. The transformation $\g$ is said to be \emph{mixing} if for any $f$,~$g \in L^2(X, \mu)$, 
\begin{equation*}
\int_X (f \circ \g^n)g \mbox{d}\mu \to \tau(f)\tau(g) \mbox{ when } |n| \to +\infty. 
\end{equation*}
\end{df}
\paragraph{\bf Example.} Let $r \in ]0, 1[$. Let $X_0 = \{0, 1\}$ with $\mu^r_0({0}) = r$ and $\mu^r_0({1}) = 1 - r$. Let $(X_r, \mu_r)$ be the probability space $\displaystyle \prod_{\Z}(X_0, \mu^r_0)$, and let $\g_r$ be the Bernoulli shift on $(X_r, \mu_r)$ defined by: 
\begin{equation*}
\forall (c_n)_{n \in \Z} \in X, \; \g_r\cdot(c_n)_{n \in \Z} = (c_{n + 1})_{n \in \Z}.
\end{equation*}
It is well-known that $\g_r$ is a measure-preserving, free, ergodic, mixing transformation. So, we obtain a measure-preserving, free, ergodic, mixing action of $\Z$ on the probability space $(X_r, \mu_r)$. Moreover, according to \cite{petersen}, the entropy of the transformation $\g_r$ is given by $H(\g_r) = -(r\log_2(r) + (1 - r)\log_2(1 - r))$. Therefore, we obtain a continuum of pairwise non-conjugate measure-preserving, free, ergodic, mixing actions.

From now on, let $(\g^n)_{n \in \Z}$ be a measure-preserving, free, ergodic, mixing action of $\Z$ on $L^{\infty}(X, \mu)$. We are going to induce this action up to $\R$ in the following way. Let 
\begin{eqnarray*}
A & = & L^{\infty}(\R \times X)^{\Z}\\
& = & \left\{F \in L^{\infty}(\R \times X), n\cdot F = F, \forall n \in \Z\right\} \\
& = & \left\{F \in L^{\infty}(\R \times X), F(t, \g^n(x)) = F(t + n, x), \forall (n, x, t) \in \Z \times X \times \R\right\}.
\end{eqnarray*}
We consider the action $(\s_t)$ of $\R$ on $A$ defined by: for any $s$, $t \in \R$, $x \in X$ and $F \in A$, $(\s_tF)(s, x) = F(s - t, x)$. It is well known according to \cite{takesaki73} that $(\s_t)$ is a measure-preserving, free, ergodic action of $\R$ on $A$, and $A \rtimes_{\s}\R$ is isomorphic to $(L^{\infty}(X) \rtimes_{\g}\Z)\otimes B(L^2(\R/\Z))$. We can identify $A$ with $L^{\infty}(Y, \nu)$ where $(Y, \nu)$ isomorphic to $([0, 1[, \la) \times (X, \mu)$ as a probability space. More precisely, for any $t \in \R$, let $t = [t] + \{t\}$, with $[t]$ the entire part of $t$. The mapping $\theta : L^2([0, 1[, \la) \otimes L^2(X, \mu) \to L^2(Y, \nu)$ defined by, $\theta(\xi \otimes \eta)(t, x) = \xi(\{t\})\eta(\g^{[t]}(x))$, for $\xi \in L^2([0, 1[, \la)$ and $\eta \in L^2(X, \mu)$ is an isomorphism of Hilbert spaces. Through this isomorphism, we can identify $L^2([0, 1[, \la) \otimes \C1$ with $L^2(\R/\Z, \la)$. We shall denote by $\h_0$ the orthogonal of $L^2(\R/\Z)$ in $L^2(Y, \nu)$, i.e. $\h_0$ is nothing but $L^2([0, 1[) \otimes (L^2(X, \mu) \ominus \C1)$. Moreover, the action $(\s_t)$ of $\R$ on $A = L^{\infty}(Y)$ gives rise to a unitary representation $(U_t)$ of $\R$ on the Hilbert space $L^2(Y, \nu)$: if the canonical embedding of $L^{\infty}(Y)$ into $L^2(Y, \nu)$ is denoted by $\eta$, $(U_t)$ is defined by $U_t\eta(y) = \eta(\s_{-t}(y))$ for any $y \in L^{\infty}(Y)$.  We want to go further and prove the following result: 
\begin{prop}\label{mixing}
With the previous notations, we have for any $\zeta_1$, $\zeta_2 \in \h_0$, $\langle\zeta_1, U_t\zeta_2\rangle \to 0$ when $|t| \to +\infty$.
\end{prop}
\begin{proof}
It suffices to show that for any $\xi_1$, $\xi_2 \in L^2([0, 1[, \la)\cap L^{\infty}([0, 1[)$ and any $\eta_1$, $\eta_2 \in (L^2(X, \mu) \ominus \C1) \cap L^{\infty}(X)$, $\langle\xi_1\otimes\eta_1, U_t(\xi_2\otimes\eta_2)\rangle \to 0$ when $|t| \to +\infty$. Let $\xi_1$, $\xi_2 \in L^2([0, 1[, \la)\cap L^{\infty}([0, 1[)$, $\eta_1$, $\eta_2 \in (L^2(X, \mu) \ominus \C1)\cap L^{\infty}(X)$ and $t \in \R$. We get
\begin{eqnarray}
\langle\xi_1\otimes\eta_1, U_t(\xi_2\otimes\eta_2)\rangle & = & \iint_{[0, 1[ \times X}\xi_1(\{s\})\eta_1(\g^{[s]}(x))\xi_2(\{s + t\})\eta_2(\g^{[s + t]}(x)) \mbox{d}\mu(x)\mbox{d}\la(s) \nonumber \\
& = & \int_{[0, 1[}\xi_1(\{s\})\xi_2(\{s + t\})\left(\int_{X}\eta_1(\g^{[s]}(x))\eta_2(\g^{[s + t]}(x)) \mbox{d}\mu(x)\right)\mbox{d}\la(s). \label{convergence}
\end{eqnarray}
As the action $(\g^{n})_{n \in \Z}$ is mixing, it is clear that with $s \in \R$ fixed,  $$\int_{X}\eta_1(\g^{[s]}(x))\eta_2(\g^{[s + t]}(x)) \mbox{d}\mu(x) \to 0$$ when $|t| \to +\infty$, because $\eta_1$, $\eta_2 \in (L^2(X, \mu) \ominus \C1)\cap L^{\infty}(X)$. The function $f : s \mapsto \xi_1(\{s\})\xi_2(\{s + t\})\left(\int_{X}\eta_1(\g^{[s]}(x))\eta_2(\g^{[s + t]}(x)) \mbox{d}\mu(x)\right)$ is such that for any $s \in \R$, $|f(s)| \leq \|\xi_1\| \|\xi_2\| \|\eta_1\| \|\eta_2\|$ which is integrable on $[0, 1[$ with respect to the Lebesgue measure $\la$. Finally, thanks to the dominated convergence theorem applied to $f$ with Equation $(\ref{convergence})$, we have proved that $\langle\xi_1\otimes\eta_1, U_t(\xi_2\otimes\eta_2)\rangle \to 0$ when $|t| \to +\infty$.
\end{proof}
For the flow $(\s_t)_{t \in \R}$ on the probability space $(Y, \nu)$, let $\tau(\s)$ be the weakest topology on $\R$ that makes the map from $\R$ to $\Aut(L^{\infty}(Y))$, which sends $t$ onto $\s_t$, continuous (notice that $\Aut(L^{\infty}(Y))$ is endowed with the $u$-topology).
We can prove the following result:
\begin{prop}
For the flow $(\s_t)_{t \in \R}$ on the probability space $(Y, \nu)$ as before, the topology $\tau(\s)$ is the usual topology on~$\R$.
\end{prop}
\begin{proof}
Let $(t_k)_{k \in \N}$ be a sequence of real numbers such that $t_k \to 0$ with respect to the topology $\tau(\s)$. Then, $\s_{t_k} \to \Id$ in $\Aut(L^{\infty}(Y))$ with respect to the $u$-topology. In particular, for any $\zeta \in L^2(Y)\ominus L^2(\R/\Z)$, $\zeta \neq 0$, $\langle\zeta, U_{t_k}\zeta\rangle \to \|\zeta\|^2 \neq 0$ when $k \to +\infty$. Consequently, we get, thanks to Proposition \ref{mixing}, that $(t_k)$ is necessarily bounded. Moreover, for any cluster point $t$ of the sequence $(t_k)$, we must have $\s_t = \Id$. As the flow $(\s_s)_{s \in \R}$ is free, $t = 0$ necessarily. Therefore, $t_k \to 0$ w.r.t the usual topology on $\R$.
\end{proof}
At last, it is easy to prove the following claim: let $\s_1$ and $\s_2$ be two flows which come from two actions of $\Z$, $\g_1$ and $\g_2$; if the flows $\s_1$ and $\s_2$ are conjugate, then the actions $\g_1$ and $\g_2$ are conjugate. Therefore, we have proved the following result:
\begin{theo}
Let $N$ be $L(\mathbf{F}_{\infty}) \otimes B(\h)$ endowed with the one-parameter automorphism group $(\a_t)_{t \in \R}$ scaling the trace from R\u{a}dulescu \cite{radulescu1991}. For $r \in ]0, 1/2]$, let $\g_r$ be the Bernoulli shift on $(X_r, \mu_r)$ defined as before, and $\s_r$ the flow on $L^{\infty}(Y_r)$ obtained from $\g_r$ after induction to $\R$. With the family $(N \indice{\Join}_{\s_r} L^{\infty}(Y_r))_{r \in ]0, 1/2]}$, we get an uncountable family of pairwise non-isomorphic factors of type ${\rm III_1}$. All these factors are non-full but they have the same $\tau$ invariant: it is the usual topology on $\R$. In particular, they have no almost periodic weights.
\end{theo}

\section{Appendix}

In this last section, we are going to prove a kind of \textquotedblleft $14\varepsilon$ lemma\textquotedblright $\,$ which can be viewed as a new version of Lemma $4.1$ in \cite{vaes2004} and Lemma $4.1$ in \cite{vaes2005}, which were also a generalization of $14\varepsilon$ lemma due to Murray $\&$ von Neumann \cite{mvn} (see \cite{barnett95} for another version of this lemma in case of a free product of two von Neumann algebras of type ${\rm III}$). We will denote by $(P, \omega)$ a free Araki-Woods factor with its free quasi-free state $\omega$: we remind that $\omega$ is given by $\left\langle \cdot \Omega, \Omega\right\rangle$ with $\Omega$ the vacuum vector (see \cite{shlya97} for further details). We shall denote by $A$ the von Neumann algebra $L^{\infty}(X) \otimes B(L^2(\R))$. Let $\s$ be the canonical embedding of $(L^{\infty}(X) \otimes P) \rtimes_{\g \otimes \s^{\omega}} \R$ into $P \otimes A$; in particular, for any $y \in P$, $\s(y) = (\s^{\omega}_{-t}(y))_{t \in \R}$ and thus $\s(y) \in P \otimes L^{\infty}(X)$. Let $E : P \otimes A \to A$ be the conditional expectation defined by $E = \omega \otimes \id$. Let $\k = L^2(X) \otimes L^2(\R)$ be the Hilbert space on which $A$ acts; let $\xi \in \k$ such that $\|\xi\| = 1$, and $\omega_{\xi}$ the vector state on $A$ associated with $\xi$. We shall denote by $\left\| \cdot \right\|_{\omega \otimes \omega_{\xi}}$ the semi-norm with respect to the normal state $\omega \otimes \omega_{\xi}$ on $P \otimes A$.

\begin{lem}\label{14epsilon}
Let $(P_i, \omega_i)$ be two von Neumann algebras endowed with a faithful normal state $\omega_i$ ($i = 1, 2$) such that we can write $(P, \omega) = (P_1, \omega_1) \ast (P_2, \omega_2)$. Let $a \in P_1$ and $b$, $c \in P_2$. Assume that $a$, $b$, $c$ are analytic with respect to the state $\omega$. For any $x \in P \otimes A$ and any $\xi \in \k$, $\|\xi\| = 1$,
\begin{eqnarray*}
\left\|x - E(x)\right\|_{\omega \otimes \omega_{\xi}}  & \leq & \mathcal{E}(a, b, c)\max\left\{\|[\s(a), x]\|_{\omega \otimes \omega_{\xi}}, \|[\s(b), x]\|_{\omega \otimes \omega_{\xi}}, \|[\s(c), x]\|_{\omega \otimes \omega_{\xi}}\right\}\\
& &  + \mathcal{F}(a, b, c)\left\|x\right\|
\end{eqnarray*}
with 
\begin{eqnarray*}
\mathcal{E}(a, b, c) & = & 6\|a\|^3 + 4\|b\|^3 + 4\|c\|^3,\\ 
\mathcal{F}(a, b, c) & = & 3\mathcal{C}(a) + 2\mathcal{C}(b) + 2\mathcal{C}(c) + 12|\omega(c^*b)|\:\|c^*b\|,\\ 
\mathcal{C}(a) & = & 2\|a\|^3\|a - \s^{\omega}_{i/2}(a)\| + 2\|a\|^2\|\s(a)(a^*\otimes1) - 1\|^2_{\omega\otimes\omega_{\xi}}\\ 
& &  + 3(1 + \|a\|^2)\|a^*a - 1\| + 6|\omega(a)|\:\|a\|.
\end{eqnarray*}
\end{lem}
\begin{proof}
For each $P_i$ let us denote by $\h_i$ the Hilbert space which comes from the GNS representation of $\omega_i$ and let $\xi_i$ be the associated cyclic vector. Let us take $(\h, \Omega) = (\h_1, \Omega_1) \ast (\h_2, \Omega_2)$. We remind \cite{voiculescu92} that
$$\h = \C\Omega \oplus (\mathring{\h}_1 \otimes \h(2, l)) \oplus (\mathring{\h}_2 \otimes \h(1, l))$$
with $\mathring{\h}_i = \h_i \ominus \C\Omega_i$, 
$$\h(2, l) = \C\Omega \oplus \mathring{\h}_2 \oplus (\mathring{\h}_2 \otimes \mathring{\h}_1) \oplus (\mathring{\h}_2 \otimes \mathring{\h}_1 \otimes \mathring{\h}_2) \oplus \cdots,$$
$$\h(1, l) = \C\Omega \oplus \mathring{\h}_1 \oplus (\mathring{\h}_1 \otimes \mathring{\h}_2) \oplus (\mathring{\h}_1 \otimes \mathring{\h}_2 \otimes \mathring{\h}_1) \oplus \cdots.$$
Moreover we shall denote by $\widetilde{\h}_1$ and $\widetilde{\h}_2$ the Hilbert spaces $\mathring{\h}_1 \otimes \h(2, l) \otimes \k$ and $\mathring{\h}_2 \otimes \h(1, l) \otimes \k$.

Let $x \in P \otimes A$ and let us define $\eta = x\cdot(\Omega \otimes \xi)$. We write $\eta = (1 \otimes E(x))\cdot(\Omega \otimes \xi) + \mu + \g$ with $\mu \in \widetilde{\h}_1$ and $\g \in \widetilde{\h}_2$. We define for $\zeta \in \h$ and $y \in P$, $\zeta\cdot y = Jy^*J\cdot\zeta$, and we notice that $(z\Omega)\cdot\s^{\omega}_{i/2}(y) = zy\Omega$ for $y \in D(\s^{\omega}_{i/2})$. Let $\mathring{x} = x - 1\otimes E(x)$, $\eta_0 = \mu + \g$, $\tilde{\eta} = \s(a)\cdot\eta\cdot(a^*\otimes1)$, $\tilde{\g} = \s(a)\cdot\g\cdot(a^*\otimes1)$, $\tilde{\mu} = \s(a)\cdot\mu\cdot(a^*\otimes1)$ and $\tilde{\zeta} = \eta_0 - \g - \tilde{\g}$. First of all, we are going to assess the quantity $\left\|\tilde{\eta} - \eta\right\|_{\omega \otimes \omega_{\xi}}$. We have
\begin{eqnarray*}
\tilde{\eta} & = & \s(a)\cdot\eta\cdot(a^*\otimes1)\\
                   & = & \s(a)x\cdot(\Omega\otimes\xi)\cdot(a^*\otimes1)\\
                   & = & [\s(a), x]\cdot(\Omega\otimes\xi)\cdot(a^*\otimes1) + x\s(a)\cdot(\Omega\otimes\xi)\cdot(a^*\otimes1)\\
                   & = & [\s(a), x]\cdot(\Omega\otimes\xi)\cdot(a^*\otimes1) + x\s(a)\cdot(\Omega\otimes\xi)\cdot((a^* - \s^{\omega}_{i/2}(a^*))\otimes1)\\ 
                     &  &  + x\s(a)\cdot(\Omega\otimes\xi)\cdot(\s^{\omega}_{i/2}(a^*)\otimes1)\\
                   & = & [\s(a), x]\cdot(\Omega\otimes\xi)\cdot(a^*\otimes1) + x\s(a)\cdot(\Omega\otimes\xi)\cdot((a^* - \s^{\omega}_{i/2}(a^*))\otimes1)\\ 
                    &  &  + x\s(a)(a^*\otimes1)\cdot(\Omega\otimes\xi).
\end{eqnarray*}
Thus, we get immediatly
\begin{eqnarray}\label{eta}
\left\|\tilde{\eta} - \eta\right\|_{\omega \otimes \omega_{\xi}} & \leq & \|a\|\left\|[\s(a), x]\right\|_{\omega \otimes \omega_{\xi}} + \|x\|\|a\|\|a - \s^{\omega}_{i/2}(a)\|\\
& & + \|x\|\left\|\s(a)(a^*\otimes1) - 1\right\|_{\omega \otimes \omega_{\xi}}.\nonumber
\end{eqnarray}
Moreover, a straightforward computation gives us 
\begin{eqnarray*}
\left| \|\g\|^2 - \|\tilde{\g}\|^2\right| & \leq & (1 + \|a\|^2)\|a^*a - 1\|\left\|\mathring{x}\right\|^2_{\omega \otimes \omega_{\xi}},\\
|\langle\tilde{\mu}, \tilde{\g}\rangle| & \leq & (1 + \|a\|^2)\|a^*a - 1\|\left\|\mathring{x}\right\|^2_{\omega \otimes \omega_{\xi}}.
\end{eqnarray*}
Let us denote by $Q_2$ the projection onto the Hilbert space $\tilde{\h}_2$. As $a \in P_1$, we can easily see that
 \begin{equation*}
\|Q_2(\s(a)(JaJ\otimes1)\cdot\g)\| \leq |\omega(a)|\|(JaJ\otimes1)\cdot\g\|.
\end{equation*}
Thus, we obtain the following inequality
\begin{equation*}
|\langle\tilde{\g}, \g\rangle| \leq |\omega(a)|\|a\|\|\mathring{x}\|^2_{\omega\otimes\omega_{\xi}}.
\end{equation*}
Because $\mu \perp \eta$, we get immediately $|\langle\tilde{\zeta}, \g\rangle| = |\langle\tilde{\g}, \g\rangle|$. We still have to assess the quantity $|\langle\tilde{\zeta}, \tilde{\g}\rangle|$. Before doing that, we can notice that as $a \in P_1$, $\s(a^*)(1\otimes E(x))(\Omega \otimes \xi) \in \C \Omega \oplus \tilde{\h}_1$. Hence
\begin{eqnarray*}
\langle\eta - \eta_0, \tilde{\g}\rangle & = & \langle(1\otimes E(x))\cdot(\Omega \otimes \xi), \s(a)\cdot\g\cdot(a^*\otimes1)\rangle\\
                                                                & = & \langle\s(a^*)(1\otimes E(x))\cdot(\Omega \otimes \xi), \g\cdot(a^*\otimes1)\rangle\\ 
                                                                & = & 0.\label{premiere}
\end{eqnarray*}
We can prove exactly in the same way that
\begin{equation*}
\langle\tilde{\eta} - \tilde{\eta_0}, \tilde{\g}\rangle = 0.
\end{equation*}
Now
\begin{eqnarray*}
\langle\tilde{\zeta}, \tilde{\g}\rangle & = & \langle\eta_0 - \g - \tilde{\g}, \tilde{\g}\rangle\\
                                                              & = & \langle\eta - \g - \tilde{\g}, \tilde{\g}\rangle\\
                                                              & = & \langle\eta - \tilde{\eta}, \tilde{\g}\rangle - \langle\g, \tilde{\g}\rangle + \langle\tilde{\eta} - \tilde{\g}, \tilde{\g}\rangle\\
                                                              & = & \langle\eta - \tilde{\eta}, \tilde{\g}\rangle - \langle\g, \tilde{\g}\rangle + \langle\tilde{\eta_0} - \tilde{\g}, \tilde{\g}\rangle\\
                                                              & = & \langle\eta - \tilde{\eta}, \tilde{\g}\rangle - \langle\g, \tilde{\g}\rangle + \langle\tilde{\mu}, \tilde{\g}\rangle
\end{eqnarray*}
because $\langle\eta - \eta_0, \tilde{\g}\rangle = \langle\tilde{\eta} - \tilde{\eta_0}, \tilde{\g}\rangle = 0$. Thus, thanks to Equation $(\ref{eta})$, we get an inequality for $|\langle\tilde{\zeta}, \tilde{\g}\rangle|$. Of course
\begin{eqnarray*}
\|\mu\|^2 + \|\g\|^2 & = & \|\tilde{\zeta} + \g + \tilde{\g}\|^2\\ 
& \geq & \|\tilde{\zeta}\|^2 + \|\g\|^2 + \|\tilde{\g}\|^2 - 2|\langle\tilde{\zeta}, \g\rangle| - 2|\langle\tilde{\zeta}, \tilde{\g}\rangle| - 2|\langle\g, \tilde{\g}\rangle|\\
                                & \geq & 2\|\g\|^2 - \left| \|\g\|^2 - \|\tilde{\g}\|^2\right| - 2|\langle\tilde{\zeta}, \g\rangle| - 2|\langle\tilde{\zeta}, \tilde{\g}\rangle| - 2|\langle\g, \tilde{\g}\rangle|.
\end{eqnarray*}
Finally, noticing that $\|\mathring{x}\|_{\omega\otimes\omega_{\xi}} \leq \|x\|_{\omega\otimes\omega_{\xi}}$ because $E = \omega \otimes \id$, and using all the previous inequalities we obtained above, we have
\begin{equation}\label{equationa}
\|\g\|^2 \leq \|\mu\|^2 + 2 \|a\|^3\|[\s(a), x]\|_{\omega\otimes\omega_{\xi}} \|\mathring{x}\|_{\omega\otimes\omega_{\xi}} + \mathcal{C}(a)\|x\| \|\mathring{x}\|_{\omega\otimes\omega_{\xi}}.
\end{equation}

We can do now exactly the same thing as before with $b$ and $c$ instead of $a$. Indeed, let $\eta' = \s(b)\cdot\eta\cdot(b^*\otimes1)$, $\eta'' = \s(c)\cdot\eta\cdot(c^*\otimes1)$, $\mu' = \s(b)\cdot\mu\cdot(b^*\otimes1)$ and $\mu'' = \s(c)\cdot\mu\cdot(c^*\otimes1)$. We define $\zeta' = \eta_0 - \mu - \mu' - \mu''$. We find that
\begin{eqnarray*}
\|\mu\|^2 + \|\g\|^2 & \geq & 3\|\mu\|^2 - \left| \|\mu\|^2 - \|\mu'\|^2\right| - \left| \|\mu\|^2 - \|\mu''\|^2\right| - 2|\langle\zeta', \mu\rangle| - 2|\langle\zeta', \mu'\rangle| \\
                                &    & - 2|\langle\zeta', \mu''\rangle| - 2|\langle\mu, \mu'\rangle| - 2|\langle\mu, \mu''\rangle| - 2|\langle\mu', \mu''\rangle|.
\end{eqnarray*}
Once again, we assess the negative terms and we get
\begin{eqnarray}\label{equationbc}
2\|\mu\|^2  & \leq &  \|\g\|^2 + 2 \|b\|^3\|[\s(b), x]\|_{\omega\otimes\omega_{\xi}} + 2 \|c\|^3\|[\s(c), x]\|_{\omega\otimes\omega_{\xi}} \|\mathring{x}\|_{\omega\otimes\omega_{\xi}} \\ 
                    &    & + \left(\mathcal{C}(b) + \mathcal{C}(c) + 6 |\omega(c^*b)|\|c^*b\|\right)\|x\| \|\mathring{x}\|_{\omega\otimes\omega_{\xi}}. \nonumber
\end{eqnarray}
As $\|\mu\|^2 + \|\g\|^2 = \|\mathring{x}\|^2_{\omega\otimes\omega_{\xi}}$, a combination of inequalities (\ref{equationa}) and (\ref{equationbc}) gives the inequality of the lemma.
\end{proof}  
Let $P_i$ be two free Araki-Woods factors endowed with their free quasi-free state $\omega_i$, $(i = 1, 2)$ such that $(P, \omega) = (P_1, \omega_1) \ast (P_2, \omega_2)$. We know, thanks to Lemma $4.3$ in \cite{vaes2004}, that $P_1$ contains a bounded sequence of elements $(a_n)$ analytic w.r.t. the state $\omega$ and which satisfy $\|\s^{\omega}_{i/2}(a_n) - a_n\| \to 0$, $\|a_n^*a_n - 1\| \to 0$, $\omega(a_n) \to 0$ and $\s(a_n)(a_n^*\otimes1) - 1 \to 0$  $\ast-$strongly.
We know besides that $P_2$ contains bounded sequences $(b_n)$ and $(c_n)$ which satisfy the same condition as $(a_n)$ and the condition $\omega(c_n^*b_n) \to 0$. We can now state the following proposition, used in proof of Lemma \ref{technique}:
\begin{prop}\label{conclusion}
Let $(x_k)$ be a bounded sequence in $C\stackrel{\s}{\subset} M = P \otimes L^{\infty}(X) \otimes B(L^2(\R))$ which almost commutes with $\s(P)$. Then $E(x_k) - x_k \to 0$ $\ast-$strongly.
\end{prop}
\begin{proof}
The proof of this proposition is a straightforward application of Lemma $\ref{14epsilon}$. Let $\xi \in \k$, $\|\xi\| = 1$,  $\varepsilon > 0$, and $M = \sup\left\{\|a_n\|, \|b_n\|, \|c_n\|, \|x_k\|\right\}$. As $\s(a_n)(a_n^*\otimes1) - 1 \to 0$ $\ast$-strongly, it is clear that $\|\s(a_n)(a_n^*\otimes1) - 1\|_{\omega\otimes\omega_{\xi}} \to 0$. So, there exists $n \in \N$ such that $\mathcal{F}(a_n, b_n, c_n) < \frac{\varepsilon}{M}$. Then with this particular $n$, as $(x_k)$ almost commutes with $\s(P)$, there exists $k_0 \in \N$ such that for any $k \geq k_0$, $$\max\left\{\left\|[\s(a_n), x_k]\right\|_{\omega \otimes \omega_{\xi}}, \left\|[\s(b_n), x_k]\right\|_{\omega \otimes \omega_{\xi}}, \left\|[\s(c_n), x_k]\right\|_{\omega \otimes \omega_{\xi}}\right\} < \frac{\varepsilon}{14M^3}.$$ Hence, we have proved that
\begin{equation*}
\forall \xi \in \k \; \forall \varepsilon > 0 \; \exists \; k_0 \in \N \; \forall k \geq k_0 \; \left\|x_k - E(x_k)\right\|_{\omega \otimes \omega_{\xi}}  \leq \varepsilon.
\end{equation*}
That means exactly that $x_k - E(x_k) \to 0$ strongly in $M$. We can do the same thing for $x_k^*$. Thus, we have proved that $x_k - E(x_k) \to 0$ $\ast-$strongly in $M$.
\end{proof}

\section*{Acknowledgment}

We wish to gratefully thank our advisor Stefaan Vaes for many helpful discussions and many encouragements during the preparation of this paper.

\addcontentsline{toc}{section}{Bibliography}
\bibliographystyle{plain}

\begin{thebibliography}{AA}


\bibitem{barnett95} {\sc L. Barnett}, {\it Free product von Neumann algebras of type ${\rm III}$}. Proc. Amer. Math. Soc. {\bf123} (1995), 543--553.

\bibitem{connes74} {\sc A. Connes}, {\it Almost periodic states and factors of type ${\rm III_1}$}. J. Funct. Anal. {\bf 16} (1974), 415--445.

\bibitem{connes76} {\sc A. Connes}, {\it Classification of injective factors.} Ann. of Math. {\bf 104} (1976), 73--115.

\bibitem{connes73} {\sc A. Connes},
{\it Une classification des facteurs de type ${\rm III}$.} Ann. Sci. {\'E}cole Norm. Sup. {\bf 6} (1973), 133--252.








\bibitem{mvn} {\sc F. Murray \& J. von Neumann}, {\it Rings of operators.} ${\rm IV}$.  Ann. of Math. {\bf 44} (1943), 716--808.

\bibitem{OP} {\sc N. Ozawa \& S. Popa}, { \it Some prime factorization results for type ${\rm II_1}$ factors.}  Invent. Math. {\bf 156} (2004), 223--234.

\bibitem{petersen} {\sc K. Petersen}, {\it Ergodic theory.} Cambridge Studies in Advanced Mathematics. {\bf 2}. Cambridge University Press, Cambridge,1983.


\bibitem{popa2001} {\sc S. Popa}, {\it On a class of type ${\rm II_1}$ factors with Betti numbers invariants. } Ann. of Math. {\bf 163} (2006), 809--899


\bibitem{popa2003} {\sc S. Popa}, {\it On the distance between MASA's in type ${\rm II_1}$ factors.} Fields Institute Communications. {\bf 2} (2001), 321--324.

\bibitem{popaduff} {\sc S. Popa}, {\it On Ozawa's property for free group factors.} Preprint $2006$. math.OA/0608451.



\bibitem{radulescu1991} {\sc F. R\u{a}dulescu}, { \it A one-parameter group of automorphisms of $L(\mathbf{F}_{\infty}) \otimes B(H)$ scaling the trace.} C. R. Acad. Sci. Paris Sér. ${\rm I}$ Math. {\bf 314} (1992), 1027--1032.

\bibitem{schmidt1980} {\sc K. Schmidt}, {\it Asymptotically invariant sequences and an action of SL$(2, \Z)$ on the $2$-sphere.}  Israel J. Math. {\bf 37} (1980), 193--208.

\bibitem{schmidt1981} {\sc K. Schmidt}, {\it Amenability, Kazhdan property T, strong ergodicity and invariant means for ergodic group actions.} Ergod. Th. \& Dynam. Sys. {\bf 1} (1981), 223--236.

\bibitem{shlya2003} {\sc D. Shlyakhtenko}, {\it On multiplicity and free
  absorption for free Araki-Woods factors.} Preprint $2003$. math.OA/0302217.

\bibitem{shlya2004} {\sc D. Shlyakhtenko}, {\it On the classification of full
  factors of type {\rm III}.} Trans. of the Amer. Math. Soc. {\bf 356} (2004), 4143--4159.


\bibitem{shlya98} {\sc D. Shlyakhtenko}, {\it Some applications of freeness with amalgamation.} J. Reine Angew. Math. {\bf 500} (1998), 191--212.

\bibitem{shlya97} {\sc D. Shlyakhtenko}, {\it Free quasi-free states.} Pacific J. Math. {\bf 177} (1997), 329--368.

\bibitem{takesakiII} {\sc M. Takesaki}, { \it Theory of Operator Algebras ${\rm II}$.} {\it EMS} {\bf 125}. Springer-Verlag, Berlin, Heidelberg, New-York, 2000.


\bibitem{takesaki73} {\sc M. Takesaki}, {\it Duality for crossed products and structure of von Neumann algebras of type ${\rm III}$.} Acta Math. {\bf 131} (1973), 249--310.

\bibitem{vaes2004} {\sc S. Vaes}, {\it \'Etats quasi-libres libres et facteurs de
  type {\rm III} (d'apr{\`e}s D. Shlyakhtenko).} S{\'e}minaire Bourbaki, expos\'e 937, {Ast\'erisque {\bf
  299}} (2005), 329--350.

\bibitem{vaes2005} {\sc S. Vaes}, {\it Strictly outer actions of groups and quantum groups.} J. reine angew. Math. {\bf  578} (2005), 147--184.

\bibitem{voiculescu92} {\sc D.V. Voiculescu, K.J. Dykema \& A. Nica}, {\it Free
  random variables.} CRM Monograph Series {\bf 1}.
American Mathematical Society, Providence, RI, 1992. 




\end{thebibliography}

\end{document}